\documentclass[aos,MSNbibl,dvips]{arximspdf}
\usepackage{mathbh}

\doi{10.1214/14-AOS1253} %kopijuoti is PTS
\volume{42}
\issue{5}
\pubyear{2014}
\firstpage{2058}
\lastpage{2091}
\docsubty{FLA}

\makeatletter
\newcommand{\rrvert}{\vert}
\newcommand{\rrVert}{\Vert}
\newcommand{\llvert}{\vert}
\newcommand{\llVert}{\Vert}
\newcommand{\given}{|}
\newtheorem{teo}{Theorem}
\newtheorem{prop}{Proposition}
\newtheorem{lem}{Lemma}
\newproclaim{remark}{Remark}
\makeatother

\begin{document}
\begin{frontmatter}

\title{On Bayesian supremum norm contraction rates}% of
%posterior distributions}
\runtitle{On Bayesian sup-norm rates}

\begin{aug}
\author{\fnms{Isma\"{e}l}~\snm{Castillo}\corref{}\ead[label=e1]{ismael.castillo@upmc.fr}\thanksref{t1}}
\runauthor{I. Castillo}
\affiliation{CNRS---LPMA Paris}
\address{CNRS---Laboratoire Probabilit\'es\\
\quad et Mod\`eles Al\'eatoires\\
%LPMA - UMR 7599\\
Universit\'es Paris VI \& VII \\
B\^{a}timent Sophie Germain\\
75205 Paris Cedex 13\\
France\\
\printead{e1}} %adresu isvedimo komanda gale!
\end{aug}
\thankstext{t1}{Supported in part by ANR Grant ``Banhdits''
ANR-2010-BLAN-0113-03.}

% HISTORY:
\received{\smonth{4} \syear{2013}}
\revised{\smonth{3} \syear{2014}}

% ABSTRACT
%
\begin{abstract}
Building on ideas from Castillo and Nickl
[\textit{Ann. Statist.} \textbf{41} (2013) \mbox{1999--2028}],
a method is
provided to study nonparametric Bayesian posterior convergence rates
when ``strong'' measures of distances, such as the sup-norm, are
considered. In particular, we show that likelihood methods can achieve
optimal minimax sup-norm rates in density estimation on the unit
interval. The introduced methodology is used to prove that commonly
used families of prior distributions on densities, namely log-density
priors and dyadic random density histograms, can indeed achieve optimal
sup-norm rates of convergence. New results are also derived in the
Gaussian white noise model as a further illustration of the presented
techniques.
\end{abstract}

% KEYWORDS
\begin{keyword}[class=AMS]
\kwd[Primary ]{62G20}
\kwd[; secondary ]{62G05}
\kwd{62G07}
\end{keyword}
\begin{keyword}
\kwd{Bayesian nonparametrics}
\kwd{contraction rates}
\kwd{supremum norm}
\end{keyword}
\end{frontmatter}

%s1 #&#
\section{Introduction}\label{sec1}

In the fundamental contributions by Ghosal, Ghosh and van der Vaart
\cite{ggv00}, Shen and Wasserman \cite{sw01} and Ghosal and van der
Vaart \cite{gvni}, a~general theory is developed to study the
behaviour of Bayesian posterior distributions. A main tool is provided
by the existence of exponentially powerful tests between a point and
the complement of a ball for some distance. The use of some important
distances, such as the Hellinger distance between probability measures,
indeed guarantees the existence of such tests. The theory often also
allows extensions to other metrics, for instance, $L^2$-type distances,
but the question of dealing with arbitrary metrics has been left
essentially open so far. Although a general theory might be harder to
obtain, it is natural to consider such a problem in simple, canonical,
statistical settings first, such as Gaussian white noise or density
estimation. This is the starting point of the authors in Gin\'e and
Nickl \cite{gn11}, and this paper was the first to provide tools to
get rates in strong norms, such as the $L^\infty$-norm. Exponential
inequalities for frequentist estimators are used in \cite{gn11} as a
way to build appropriate tests, and this enables one to obtain some
rates in sup-norm in density estimation. In the case where the true
density is itself supersmooth and a kernel mixture is used as a prior,
the nearly parametric minimax rate is attained, at least up to a
possible logarithmic term; see also the work by Scricciolo \cite{s12} for
related results. In the general case where the true density belongs to
a H\"{o}lder class, a sup-norm rate is obtained which differs from the
minimax rate by a power of $n$.
%, which are however slower than the minimax rates.
On the other hand, by using explicit computations, the authors in \cite
{gn11} show that in the Gaussian white noise model with conjugate
Gaussian priors, minimax sup-norm rates are attainable, which leads to
the natural question to know whether this is still possible in density
estimation, or in nonconjugate regression settings. This nontrivial
question also arises for other likelihood methods, such as
nonparametric maximum likelihood estimation; see Nickl \cite{n07}.

From a general statistical perspective, density estimation in
supremum-norm is a central problem both from theoretical and practical
points of view. The problem was the object of much interest in the
framework of minimax theory. Lower bounds in density estimation in
sup-norm can be found in Hasminskii \cite{hasminskii78}, upper-bounds
in Ibragimov and Hasminskii \cite{ih80} for density estimation and
Stone \cite{stone82} for regression. %Recent advances have pushed
%these properties even further, to the level of exact constants, see
%e.g. \cite{korostelev93} and \cite{tsy98} and to more general classes,
We refer to Goldenshluger and Lepski \cite{goldlep13} for an overview
of current work in this area. From the practical perspective, sup-norm
properties are of course very desirable, since saying that two curves
in a simulation picture look close is very naturally, and often
implicitly, done in a sup-norm sense.

Here, we establish that minimax optimal sup-norm rates of convergence
in density estimation are attainable by common and natural Bayes
procedures.
The methodology we introduce is in fact related to a programme initiated in \cite{CN12} and
continued in \cite{CN14}, namely nonparametric Bernstein--von Mises type results, as \mbox{discussed} below.
In \cite{CN14}, we use the results of the present paper to derive nonparametric Bernstein--von
Mises theorems in density estimation, as well as Donsker-type results for the posterior distribution function.
The testing approach commonly used to establish posterior rates is replaced here by tools from semiparametric Bernstein--von
Mises results (testing is still typically useful to establish
preliminary rates); see \cite{CN12} for an overview of references. We
split the distance of interest in simpler pieces, each simpler piece
being a semiparametric functional to study. One novelty of the paper
consists in providing well-chosen uniform approximation
schemes of various influence functions appearing at the semiparametric
level when estimating those simple functionals.

Two natural families of nonparametric priors are considered for
density estimation: priors on log-densities; see, for example, Ghosal, Ghosh and van der Vaart~\cite{ggv00},
Scricciolo \cite{s06}, Tokdar and Ghosh \cite{tg07}, van der Vaart
and van Zanten \cite{vvvz,ic08}, Rivoirard and Rousseau \cite{rr12},
and random (dyadic) histogram priors; see,
for example, Barron \cite{barron88}, Barron, Schervish and Wasserman
\cite{bsw99}, Walker \cite{wal04}, Ghosal and van der Vaart \cite
{gvni}, Scricciolo \cite{s07}, Gin\'e and Nickl \cite{gn11} and the
recent semiparametric treatment in \cite{cr13}. Both classes are
relevant for applications and priors of these types have been studied
from the implementation perspective; see, for example, Lenk \cite
{Lenk88}, Tokdar \cite{Surya} and references therein for the use of
logistic Gaussian process priors, and Leonard~\cite{Leonard73},
Gasparini \cite{gasparini96} for random histogram priors.
%Advantage of Gaussian priors: `easy' to simulate (a lot if algorithms
%available, even for log-priors e.g. Adams, cite Surya ?),
%semiparametrically, leads to simple and compact
%bias-control conditions, though those can be sometimes slightly
%stronger than for heavy tails, due to fast concentration of Gaussian
%tails.
% Cite also paper by Richard on MLEs as motivation for log-densities,
%if it's the type of prior he% was considering ?

New results are also derived in the Gaussian white noise model, in the
spirit of~\cite{CN12}, for nonconjugate priors.

%Let us mention that in density estimation, when the true density is
%itself a convolution and the minimax rate is then nearly parametric,
%sup-norm rates have been recently obtained in \cite{s12}.
While working on this paper, we learned from the work by Marc Hoffmann,
Judith Rousseau and Johannes Schmidt-Hieber \cite{hrs13}, which
independently obtains sup-norm properties for different priors. Their
method is different from ours, and both approaches shed light on
different specific aspects of the problem. In Gaussian white noise,
adaptive results over H\"{o}lder classes
are obtained in \cite{hrs13} for a class of sparse priors. %, which
%can be seen as Bayesian analogs of thresholding, see \cite{cv12}.
In Theorem~\ref{teo-gwn} below, the sup-norm minimax rate for fixed
regularity is obtained for canonical priors
without sparsity enforcement. The authors also give insight on the
interplay between loss function and posterior rate, as well as an upper
bound result for fairly abstract sieve-type priors, which are shown to
attain the adaptive sup-norm rate in density estimation. This is an
interesting existence result, but no method is provided to investigate
sup-norm rates for general given priors. Although for simplicity we
limit ourselves here to the fixed regularity case, the present paper
suggests such a method and demonstrates its applicability by dealing
with several commonly used classes of prior distributions. Clearly,
there is still much to do in the understanding of posterior rates for
strong measures of loss, and we hope that future contributions will go
further in the different directions suggested by both the present paper
and \cite{hrs13}.
%and the present paper and \cite{hrs13} can be seen as promising first
%steps towards a better understanding of these questions.
%[Introductory part: here]

%Let us introduce some notation used throughout the paper.
Let $L^2[0,1]$ and $L^{\infty}[0,1]$, respectively, denote the space of
square integrable functions with respect to Lebesgue measure on $[0,1]$
and the space of measurable bounded functions on $[0,1]$. Theses spaces
are equipped with their usual norms, respectively, denoted $\llVert
\cdot
\rrVert_2$
(denote by ${\langle}\cdot,\cdot{\rangle}_2$ the associated inner
product) and \mbox{$\llVert\cdot\rrVert_{\infty}$.}
Let $ {\mathcal{C}}^\alpha:= {\mathcal{C}}^{\alpha}[0,1]$
denote the class of H\"{o}lder functions on $[0,1]$
with H\"{o}lder exponent $\alpha>0$.

For any $\alpha>0$ and any $n\ge1$, denote by $\bar\varepsilon
_{n,\alpha}$ the rate
% It is the standard minimax rate over H{\chr"F6}lder-classes for the
%Hellinger-loss in density estimation on $[0,1]$.
%
%e1 #&#
\begin{equation}
\label{strate} \bar\varepsilon_{n,\alpha}:= n^{-\alpha/(2\alpha +1)}.
\end{equation}
The typical minimax rate over a ball of the H\"{o}lder space $
{\mathcal{C}}^\alpha[0,1]$, $\alpha>0$,
for the sup-norm is
%
%e2 #&#
\begin{equation}
\label{mmrate} \varepsilon_{n,\alpha}^*:= \biggl( \frac{ \log{n}}{n}
\biggr)^{\alpha/(2\alpha+1) }.
\end{equation}
Let us also set, omitting the dependence in $\alpha$ in the notation,
%
%e3 #&#
\begin{equation}
\label{def-jn} h_n = \biggl(\frac{ n }{ \log{n} } \biggr)^{-1/(2\alpha
+1)},
\qquad L_n= \bigl\lfloor\log_2 (1/h_n) \bigr
\rfloor. %\qquad\zeta_n=J_n \veps_{n,\al}.
\end{equation}

For a statistical model $\{P_f^{(n)}\}$ indexed by $f$ in some class of
functions to be specified and associated observations $X^{(n)}$, denote
by $f_0$ the ``true'' function and by $E_{f_0}^n$ the expectation under
$P_{f_0}^{(n)}$. Given a prior $\Pi$ on a set of possible
$f$'s, denote
by $\Pi[ \cdot\given X^{(n)}]$ the posterior distribution and by
$E^\Pi[ \cdot\given X^{(n)}]$ the expectation operator under the law
$\Pi[ \cdot\given X^{(n)}]$.

%s2 #&#
\section{Prologue} %- From BvM to sup-norm rates}

Let us start by a simple example in Gaussian white noise which will
serve as a
%-perhaps overly simplistic
slightly naive yet useful illustration of the main technique of proof.

Let $f$ be an element of $L^\infty[0,1]$. Let $n\ge1$. Suppose one
observes %a realisation
%
%e4 #&#
\begin{equation}
\label{gwn} dX^{(n)}(t) = f(t)\,dt + \frac{1}{\sqrt{n}}\,dW(t),\qquad t
\in[0,1],
\end{equation}
where $W$ is standard Brownian motion. Let $\{ \psi_{lk}, l\ge0, 0\le
k\le2^l-1 \}$ be a wavelet basis on the interval $[0,1]$. Here,
we take the basis constructed in \cite{CDV93}, see below for precise
definitions.
The model (\ref{gwn}) is statistically equivalent
to observing the projected observations onto the basis $\{\psi_{lk}\}$,
\[
x_{lk} = f_{lk} + \frac{1}{\sqrt{n}}\varepsilon_{lk},
\qquad l\ge0, 0\le k\le2^l-1,
\]
where $f_{lk}:={\langle}f, \psi_{lk} {\rangle}_2$ and $\varepsilon
_{lk}$ are i.i.d. standard normal.
Denote $\hat{f}_{lk}:=x_{lk}$, an efficient frequentist estimator of
the wavelet coefficient $f_{lk}$.

%s2.1 #&#
\subsection{A first example} \label{firstex}

Suppose the coefficients of the true function $f_0$ satisfy, for some
$R>0$ that we suppose to be \emph{known} in this first example,
%
%e5 #&#
\begin{equation}
\label{tolk} \sup_{l\ge0, 0 \leq k\leq2^l-1} 2^{l(1/2+\alpha)}
\llvert
f_{0,lk}\rrvert\le R, \qquad\alpha>0.
\end{equation}
Define a prior $\Pi$ on $f$ via an independent product prior on its
coordinates $f_{lk}$ onto the considered basis. The component $f_{lk}$
is assumed
to be sampled from a prior with density $\sigma_{l}^{-1}\varphi(\cdot
/\sigma_{l})$ with respect to Lebesgue measure on $[0,1]$, where, for
$\alpha, R$ as in (\ref{tolk}), $x\in\mathbb{R}$ and a given $B>R$
%
%e6 #&#
\begin{equation}
\label{prior-unif} \varphi(x) = \frac{1}{2B} \mathbh{1}_{[-B,B]}(x),
\qquad\sigma_{l}=2^{-l(1/2+\alpha)}.
\end{equation}
This type of prior was considered in \cite{gn11}, Section~2.2, and
provides a simple example of a random function with bounded $\alpha
$-H\"{o}lder norm.
%We have obtained the following proposition
%
%pr1 #&#
\begin{prop} \label{prop-unif}
Consider observations $X^{(n)}$ from the model (\ref{gwn}).
Let $f_0$~and~$\alpha$ satisfy (\ref{tolk})
and let the prior be chosen according to (\ref{prior-unif}). Then
there exists $M>0$ such that for $\varepsilon_{n,\alpha}^*$ defined
by (\ref{mmrate}),
\[
E_{f_0}^n \int\llVert f-f_0\rrVert
_{\infty} \,d\Pi\bigl(f\given X^{(n)}\bigr) \le M
\varepsilon_{n,\alpha}^*.
\]
\end{prop}
Uniform wavelet priors thus lead to the minimax rate of convergence in
sup-norm. The result has a fairly simple proof, as we now illustrate,
and is new, to the best of our knowledge.

Let $L_n$ be defined in (\ref{def-jn}). Denote by $f^{L_n}$ the
orthogonal projection of $f$
in $L^2[0,1]$ onto $\operatorname{Vect}\{\psi_{lk}, l\le L_n, 0\le
k<2^l\}$, and $f^{L_n^c}$ the\vspace*{1pt} projection of $f$ onto
$\operatorname{Vect}\{\psi_{lk}, l> L_n, 0\le k<2^l\}$. Then
\[
f-f_0 = f^{L_n} - \hat{f}{}^{L_n} +
\hat{f}{}^{L_n} - f_0^{L_n} + f^{L_n^c}
-f_0^{L_n^c},
\]
where\vspace*{1pt} $\hat{f}{}^{L_n}$ is the projection estimator onto the basis $\{
\psi_{lk}\}$ with cut-off $L_n$. Note that the previous equality as
such is an equality in $L^2$. However, if the wavelet series of $f$
into the basis $\{\psi_{lk}\}$ is absolutely convergent $\Pi$-almost
surely (which is the case for all priors considered in this paper), we
also have $f(x)=f^{L_n}(x)+f^{L_n^c}(x)$ pointwise for Lebesgue-almost
every $x$, $\Pi$-almost surely, and similarly for~$f_0$. Now,
%Decomposing as above,
%
\begin{eqnarray*}
\hspace*{-4pt}&& E^\Pi\bigl[ \llVert f-f_0\rrVert_{\infty}
\given X^{(n)} \bigr]
\\
\hspace*{-4pt}&&\qquad = \int\llVert f-f_0\rrVert
_{\infty
} \,d\Pi\bigl(f\given X^{(n)}\bigr)
\\
\hspace*{-4pt}&&\qquad \le \underbrace{\int\bigl\llVert f^{L_n} -\hat{f}{}^{L_n} \bigr
\rrVert_{\infty} \,d\Pi\bigl(f\given X^{(n)}\bigr)}_{(\mathrm{i})}
+ \underbrace{\int\bigl\llVert f^{L_n^c}\bigr\rrVert_{\infty} \,d\Pi
\bigl(f\given X^{(n)}\bigr)}_{(\mathrm{ii})} + \underbrace{\bigl\llVert
\hat{f}{}^{L_n}-f_0\bigr\rrVert_{\infty}}_{(\mathrm{iii})}.
\end{eqnarray*}
We\vspace*{1pt} have $\mbox{(iii)} \le\llVert f_0^{L_n^c}\rrVert_{\infty}+ \llVert\hat
{f}^{L_n} -
f_0^{L_n}\rrVert_{\infty}$. Using (\ref{tolk}) and the localisation
property of the wavelet basis $\llVert\sum_k \llvert\psi_{lk}\rrvert
\rrVert_\infty
\lesssim2^{l/2}$ (see below), one obtains
\[
\bigl\llVert f_0^{L_n^c}\bigr\rrVert_{\infty} \le
\sum_{l>L_n} \biggl[ \max_{k}
\llvert f_{0,lk}\rrvert\biggl\llVert\sum_{k}
\llvert\psi_{lk}\rrvert\biggr\rrVert_{\infty} \biggr]
\lesssim h_n^{\alpha} \lesssim\varepsilon_{n,\alpha}^*,
\]
where $\lesssim$ means less or equal to up to some universal constant.
The term $\llVert\hat{f}{}^{L_n} - f_0^{L_n}\rrVert_{\infty}$ depends
on the
randomness of the observations only,
\[
\bigl\llVert\hat{f}{}^{L_n} - f_0^{L_n}\bigr
\rrVert_{\infty} = \frac{1}{\sqrt{n}} \biggl\llVert\sum
_{l\le L_n, k} \varepsilon_{lk} \psi_{lk}(
\cdot)\biggr\rrVert_{\infty}.
\]
This is bounded under $E_{f_0}^n$ by a constant times $\varepsilon
_{n,\alpha}^*$; see
Lemma~\ref{lem-Gamma} for a proof in the more difficult case of
empirical processes.
\begin{longlist}[\textit{Term} (ii)]
\item[\textit{Term} (i).] By definition, $\hat f^{L_n}$ has coordinates $\hat
f_{lk}$ in the basis $\{\psi_{lk}\}$, so using the localisation
property of the wavelet basis as above, one obtains
\[
\bigl\llVert f^{L_n} - \hat{f}{}^{L_n} \bigr\rrVert
_{\infty} \lesssim\frac{1}{\sqrt{n}} \sum_{l\le L_n}
2^{l/2} \Bigl[ \max_{0\le k < 2^l} \sqrt{n}\llvert
f_{lk} - x_{lk}\rrvert\Bigr].
\]
For $t>0$, via Jensen's inequality and bounding the maximum by the sum,
using $\Pi_n$ as a shorthand notation for the posterior $\Pi[\cdot
\given X^{(n)}]$,
\begin{eqnarray*}
&& t E_{f_0}^n E^{\Pi_n}\Bigl[\max
_{0\le k < 2^l} \sqrt{n}\llvert f_{lk} - x_{lk}
\rrvert\Bigr]
\\
&&\qquad \le\log\sum_{k=0}^{2^l-1}
E_{f_0}^n E^{\Pi_n} \bigl[e^{t \sqrt
{n}(f_{lk} - x_{lk})} +
e^{-t \sqrt{n}(f_{lk} - x_{lk})} \bigr],
\end{eqnarray*}
for any $l\ge0$.
Simple computations presented in Lemma~\ref{lem-ratio} yield a
sub-Gaussian behaviour for the Laplace transform of $\sqrt
{n}(f_{lk}-x_{lk})$ under the posterior distribution, which is bounded
above by $C e^{t^2/2}$ for a constant $C$ independent of $l\le L_n$ and $k$.
From this deduce, for any $t>0$ and $l\le L_n$,
\[
E_{f_0}^n E^{\Pi}\Bigl[\max_{0\le k< 2^l}
\sqrt{n}\llvert f_{lk} - x_{lk}\rrvert \big\given
X^{(n)}\Bigr] \lesssim\frac{\log(C2^l)}{t} + \frac{t} 2.
\]
The choice $t=\sqrt{2\log(C2^l)}$ leads us to the bound
%control the expectation of the term (i) above by
%
\begin{eqnarray*}
% (i)=E_{f_0}^n \int\left\Vert f^{J_n} -f_0^{J_n} \right\Vert_{\infty}
%d\Pi(f\given
%X^{(n)})
E_{f_0}^n \mbox{(i)} & \lesssim&
\frac{1}{\sqrt{n}}\sum_{l\le L_n} \sqrt{l} 2^{l/2}
\lesssim\sqrt{L_n/(nh_n)} \lesssim\varepsilon_{n,\alpha}^*.
\end{eqnarray*}
\item[\textit{Term} (ii).] Under the considered prior, the wavelet coefficients
of $f$ are bounded by $\sigma_{l}$, so using again the localisation
property of the wavelet basis, % the term (ii) above is bounded in
%expectation by
%
\begin{eqnarray*}
E_{f_0}^n \mbox{(ii)} & \lesssim&\sum
_{l >L_n} 2^{l/2} E_{f_0}^n
E^{\Pi} \Bigl[ \max_{k} \llvert f_{lk}
\rrvert\big\given X^{(n)} \Bigr]
\\
& \lesssim&\sum_{l >L_n} 2^{l/2}
\sigma_l \lesssim h_n^{\alpha} =
\varepsilon_{n,\alpha}^*.
\end{eqnarray*}
This concludes the proof of Proposition~\ref{prop-unif}.

%%%%%%%%%%%%%%%%%%%%%%%%%%%%%%%%%%%%%%%%%%%%%%%

Although fairly simple, the previous example is revealing of some
important facts, some of which are well known from frequentist analysis
of the problem, some being specific to the Bayesian approach. The
previous proof shows two regimes of frequencies: $l\le L_n$ ``low
frequency'' and $l>L_n$ ``high frequency.'' In the low frequency
regime, the estimator $x_{lk}$ of $f_{lk}={\langle}f,\psi
_{lk}{\rangle}_2$ is satisfactory, and the concentration of the
posterior distribution around this efficient frequentist estimator is
desirable. This is reminiscent of the Bernstein--von Mises (BvM)
property; see van der Vaart
\cite{aad98}, Chapter~10, which states that in regular parametric
problems with unknown parameter $\theta$, the posterior distribution
is asymptotically Gaussian concentrating at rate $1/\sqrt{n}$ and
centered around an efficient estimator of~$\theta$.

Here are a few words on the general philosophy of the results
specifically in
the Bayesian context. Such method was used as a building block in \cite
{CN12}. The idea is to split
the distance of interest into small pieces. For the sup-norm, those
pieces can, for instance, involve the
wavelet coefficients ${\langle}f,\psi_{lk}{\rangle}_2$, but not
necessarily, as will be seen for log-density
priors. In this case, this split is obtained, for instance, from the inequality
\[
\llVert f-f_0\rrVert_{\infty} \lesssim\sum
_{l\ge0} 2^{l/2} \max_{0\le k\le
2^l -1} \bigl
\llvert{\langle}f,\psi_{lk}{\rangle}_2 - {
\langle}f_0,\psi_{lk}{\rangle}_2 \bigr
\rrvert,
\]
which holds for localised bases $\{\psi_{lk}\}$.
Note that $f\to{\langle}f,\psi_{lk}{\rangle}_2$ can be seen as a
semiparametric functional; see, for example, \cite{aad98}, Chapter~25
for an \hyperref[sec1]{Introduction} to semiparametrics and the notions of efficiency
and efficient influence functions.
Next, one analyses each piece separately, with different regimes of
indexes $l,k$ often arising, requiring specific techniques for each of them.
\begin{itemize}
\item the \emph{BvM-regime: semiparametric bias.} For ``low
frequencies,'' what is typically needed is a
concentration of the posterior distribution for the functional of
interest, say ${\langle}f,\psi_{lk} {\rangle}_2$, at rate $1/\sqrt
{n}$ around a semiparametrically efficient estimator of the
functional. This is at the
heart of the proof of semiparametric BvM results, hence the use of BvM
techniques. In particular, sharp control of the bias will be essential.
Regarding the BvM property, although the
precise Gaussian shape will not be needed here, one needs \emph{uniformity} in all frequencies in the considered regime. This requires
nontrivial strengthenings of BvM-type results, the semiparametric
efficient influence function of the functional of interest, which can
be, for instance, a re-centered version of $\psi_{lk}$, being typically
\emph{unbounded} as~$l$ grows.
\item Taking care of uniformity issues in approximation of the
efficient influence functions by the prior may require various
approximations regimes depending on $l$. For log-density priors, we
will indeed see various regimes of indexes ``$l$'' arise in the
obtained bounds for the bias.
\item The \emph{high-frequency bias} corresponds to frequencies where
the prior should make the likelihood negligible. This part can be
difficult to handle, too, especially for unbounded priors.
%One often handles this part by a direct method, but it would be
%desirable to develop more general techniques for doing this.
\end{itemize}
In the example above for uniform priors in white noise, most of the
previous steps are either almost trivial or at least can be carried out
by considering the explicit expression of the posterior, but for
different priors or in different sampling situations some of the
previous steps may become significantly harder, as we will see below.
\end{longlist}

%s2.2 #&#
\subsection{Wavelet basis and Besov spaces} \label{wave}

Central to our investigations is the tool provided by localised bases
of $L^2[0,1]$. We refer to the Lecture Notes by H\"{a}rdle,
Kerkyacharian, Picard and Tsybakov \cite{hkpt} for an \hyperref[sec1]{Introduction} to
wavelets. Two bases will be used in the sequel.

The Haar basis on $[0,1]$ is defined by $\varphi^H(x)=1$, $\psi
^H(x):=\psi^H_{0,0}(x)=-\mathbh{1}_{[0,1/2]}(x)+\mathbh
{1}_{(1/2,1]}(x)$ and
$\psi_{l,k}^H(x)=2^{l/2} \psi(2^lx-k)$, for any integer $l$ and $0\le
k\le2^l -1$. The supports of Haar wavelets form dyadic partitions of
$[0,1]$, corresponding to intervals $I_k^l:=(k2^{-l},(k+1)2^{-l}]$ for
$k>0$, and where the interval is closed to the left when $k=0$.

The boundary corrected basis of Cohen, Daubechies, Vial \cite{CDV93}
will be referred to as CDV basis. Similar to the Haar basis, the CDV
basis enables a treatment on compact intervals, but at the same time
can be chosen sufficiently smooth. A few properties are lost,
essentially simple explicit expressions, but most convenient
localisation properties and characterisation of spaces are maintained.
Below we recall some useful properties of the CDV basis. We denote this
basis $\{\psi_{lk}\}$, with indexes $l\ge0$, $0\le k\le2^l -1$ (with
respect to the original construction in \cite{CDV93}, one starts at a
sufficiently large level $l\ge J$, with $J$ fixed large enough; for
simplicity, up to renumbering, one can start the indexing at $l=0$).
Let $\alpha>0$ be fixed.
\begin{itemize}
\item$\{ \psi_{lk} \}$ forms an orthonormal basis of $L^2[0,1]$.
\item$\psi_{lk}$ have
support $S_{lk}$, with diameter at most a constant (independent of
$l,k$) times $2^{-l}$, and
$\llVert\psi_{lk}\rrVert_{\infty} \lesssim2^{l/2}$. The $\psi
_{lk}$'s are
in the H\"{o}lder class $ {\mathcal{C}}^S[0,1]$, for some $S\ge
\alpha$.
\item
At fixed level $l$, given a fixed $\psi_{lk}$ with support $S_{lk}$,
\begin{itemize}
\item[$\diamond$] the number of wavelets of the level $l'\le l$ with
support intersecting $S_{lk}$ is bounded by a universal constant
(independent of $l',l,k$),
\item[$\diamond$] the number of wavelets of the level $l'>l$ with
support intersecting $S_{lk}$ is bounded by $2^{l'-l}$ times a
universal constant. %(independent of $l',l,k$).
\end{itemize}
The following localisation property holds
$\sum_{k=0}^{2^l-1} \llVert\psi_{lk}\rrVert_{\infty} \lesssim
2^{l/2}$, where
the inequality is up to a fixed universal constant.
\item
The constant function equal to $1$ on $[0,1]$ is orthogonal to
high-level wavelets, in the sense that
${\langle}\psi_{lk},1 {\rangle}_2=\int_0^1 \psi_{lk}=0$ whenever
$l\ge M$, for a large enough constant $M$.
\item The basis $\{\psi_{lk}\}$ characterises Besov spaces
$B^{s}_{\infty,\infty}[0,1]$, any $s \le\alpha$, in terms of
wavelet coefficients. That is, $g\in B_{\infty,\infty}^s[0,1]$ if and
only if
%
%e7 #&#
\begin{equation}
\label{besovnorm} \llVert g\rrVert_{\infty,\infty,s}:= \sup_{l\ge0,
0\le k\le2^{l}-1}
2^{l(1/2+s)} \bigl\llvert{\langle}g,\psi_{lk} {
\rangle}_2\bigr\rrvert<\infty.
\end{equation}
\end{itemize}
We note that orthonormality of the basis is not essential. Other
nonorthonormal, multi-resolution dictionaries could be used instead up
to some adaptation of the proofs, as long as coefficients in the
expansion of $f$ can be recovered from inner products.
Also, recall that $B_{\infty,\infty}^s$ coincides with the H\"older
space $ {\mathcal{C}}^s$ when $s$ is not an integer and that when
$s$ is an integer the inclusion
$ {\mathcal{C}}^s\subset B_{\infty,\infty}^s$ holds. If the
Haar-wavelet is considered, the fact that $f_0$ is in $ {\mathcal
{C}}^s$, $0<s\le1$, implies that the supremum in (\ref{besovnorm})
with $\psi_{lk}=\psi_{lk}^H$ is finite.

%s3 #&#
\section{Main results}

%s3.1 #&#
\subsection{Gaussian white noise}

Consider priors $\Pi$ defined as coordinate-wise products of priors on
coordinates specified by a density $\varphi$ and scalings $\{\sigma
_l\}$ as in Section~\ref{firstex}. The next result allows for a much
broader class of priors.

Let $\varphi$ be a continuous density with respect to Lebesgue measure
on $\mathbb{R}$. We assume that
$\varphi$ is (strictly) positive on $[-1,1]$ and that it satisfies
%
%e8 #&#
\begin{equation}
\label{cprior2} \qquad\exists b_1, b_2, c_1,
c_2, \delta>0,  \forall x\dvtx  \llvert x\rrvert\ge1,\qquad
c_1 e^{-b_1\llvert x\rrvert^{1+\delta}}\le\varphi(x)\le
c_2e^{-b_2\llvert x\rrvert^{1+\delta}}.
\end{equation}
Consider a scaling $\sigma_l$ for the prior equal to, for $\delta$
the constant in (\ref{cprior2}),
%
%e9 #&#
\begin{equation}
\label{silch} \sigma_{l}= \frac{2^{-l(1/2+\alpha)}}{(l+1)^{\mu
}},\qquad\mu=
\frac{1}{1+\delta}.
\end{equation}

%th1 #&#
\begin{teo} \label{teo-gwn}
Let $X^{(n)}$ be observations from (\ref{gwn}).
Suppose $f_0$ belongs to $B_{\infty,\infty}^\alpha[0,1]$, for some
$\alpha>0$. Let the prior $\Pi$ be a product prior defined through
$\varphi$ and $\sigma_{l}$ satisfying (\ref{cprior2}), (\ref{silch}).
Then there exists $M>0$ such that for $\varepsilon_{n,\alpha}^*$
defined by (\ref{mmrate}),
\[
E_{f_0}^n \int\llVert f-f_0\rrVert
_{\infty} \,d\Pi\bigl(f\given X^{(n)}\bigr) \le M
\varepsilon_{n,\alpha}^*.
\]
\end{teo}

Theorem~\ref{teo-gwn} can be seen as a generalisation to nonconjugate
priors of Theorem 1 in \cite{gn11}. Possible choices for $\varphi$
cover several commonly used classes of prior distributions, such as
so-called exponential power (EP) distributions; see, for example, Choy
and Smith \cite{choysmith97}, Walker and Guti\'errez-Pe\~{n}a \cite
{wg99} and references therein, as well as some of the univariate
Kotz-type distibutions, see, for example, Nadarajah \cite{nad03}. %,
%but also all sub-Gaussian distributions, which includes Proposition 1
%for uniform priors as a particular case.
Other choices of prior distributions are possible, up sometimes to some
adaptations. For instance, Proposition~\ref{prop-unif} provides a
result in the case of a uniform distribution. If one allows for some
extra logarithmic term in the rate, Laplace (double-exponential)
distributions can be used, as well as distributions without the control
from below on the tail in (\ref{cprior2}), provided one chooses
$\sigma_{l}=2^{-l(1/2+\alpha)}$, as can be checked following the
steps of the proof of Theorem
\ref{teo-gwn}. As a special case, the latter include all sub-Gaussian
distributions. % can be dealt with using similar arguments up to an
%extra logarithmic factor in the rate).
Also note that Theorem~\ref{teo-gwn} as such applies to canonical
priors, in that they do not depend on $n$. Results for truncated
priors, which set ${\langle}f,\psi_{lk}{\rangle}=0$ for $l$ above a
threshold, can be obtained along the same lines, with slightly simpler proofs.

Further consequences of Theorem~\ref{teo-gwn} include the minimaxity
in sup-norm of several Bayesian estimators. The result for the
posterior mean immediately follows from a convexity argument. One can
also check that the posterior coordinate-wise median is minimax.
Details are omitted.
%are mostly for simplicity of presentation, but similar results can be
%obtained for a variety of related priors, e.g. uniform priors, Laplace
%priors (up to a logarithmic factor in the rate) etc. The reason of
%some tail control in \eqref{cprior2} is to handle priors without
%high-frequency cut-off. This way, the results are obtained for
%canonical priors, in that they do not depend on $n$. Results for
%priors with cut-off $l\le K_n$, for some $K_n\to\infty$, can be
%obtained along the same lines.

%{\it A\,daptation.} Consider the following prior $\Pi$ constructed by
%first sampling an integer $L\ge1$ according to the law $\pi(L)\propto
%e^{-L\log L}$ on the integers and, given $L$, setting
%A-Step 1. Show that for $M$ large enough,
%
%Update / Complete this according to new scheme

%s3.2 #&#
\subsection{Density estimation}

Consider independent and identically distributed observations
%
%e10 #&#
\begin{equation}
\label{mod-density} X^{(n)}=(X_1,\ldots,X_n),
\end{equation}
with unknown density function $f$ on $[0,1]$. We use the same notation
$X^{(n)}$ for observations as in the white noise model: it will always
be clear from the context which model we are referring to. Let $
{\mathcal{F}}$ be the set of densities $f$ on $[0,1]$ which are
bounded away from $0$ and $\infty$. In other words, one can write $
{\mathcal{F}}=\bigcup_{0<\rho\le D<\infty} {\mathcal{F}}(\rho,D)$,
with $ {\mathcal{F}}_{\rho,D}=\{f, 0<\rho\le f\le D<\infty, \int
_0^1f=1\}$.\vspace*{1.5pt} In the sequel, we assume that the ``true'' $f_0$
belongs to $ {\mathcal{F}}_0:= {\mathcal{F}}(\rho_0,D_0)$, for
some $0<\rho_0\le D_0<\infty$. The assumption that the density is
bounded away from $0$~and~$\infty$ is for simplicity. Allowing the
density to tend to $0$, for example, at the boundary of $[0,1]$ would
be an interesting extension, but would presumably induce technicalities
not related to our point here. %there exist constants $\rho,D$ with $0<
Let $h$ denote the Hellinger distance between densities on $[0,1]$.

%s3.2.1 #&#
\subsubsection{Log-densities priors}

Define the prior $\Pi$ on densities as follows. Given a sufficiently
smooth CDV-wavelet basis $\{\psi_{lk}\}$, consider the prior
induced by, for any $x\in[0,1]$ and $L_n$ defined in (\ref{def-jn}),
%
%e11 #&#
%e12 #&#
\begin{eqnarray}
T(x) & =& \sum_{l=0}^{L_n} \sum
_{k=0}^{2^l -1} \sigma_{lk}\alpha
_{lk} \psi_{lk}(x), \label{def-T}
\\
f(x) & =& \exp\bigl\{ T(x) - c(T) \bigr\},\qquad c(T)=\log\int
_0^1 e^{T(x)} \,dx, \label{def-prior-log}
\end{eqnarray}
where $\alpha_{lk}$ are i.i.d. random variables of density $\varphi$
with respect to Lebesgue measure on $\mathbb{R}$ and $\sigma_{lk}$
are positive reals which for simplicity we make only depend on $l$,
that is $\sigma_{lk}\equiv\sigma_l$. We consider the choices
$\varphi(x)=\varphi_G(x)=e^{-x^2/2}/\sqrt{2\pi}$ the Gaussian
density and
$\varphi(x)=\varphi_H(x)$, where $\varphi_H$ is any density such
that its logarithm $\log\varphi_H$ is Lipschitz on $\mathbb{R}$. We
refer to this as the ``log-Lipschitz case.'' For instance, the $\alpha
_{lk}$'s can be Laplace-distributed or have heavier tails, such as, for
a given $0\le\tau<1$ and
$x\in\mathbb{R}$, and $c_\tau$ a normalising constant,
%
%e13 #&#
\begin{equation}
\varphi_{H,\tau}(x) = c_\tau\exp\bigl\{-\bigl(1+\llvert x
\rrvert\bigr)^{1-\tau}\bigr\}. \label{phi-ht}
\end{equation}
Suppose the prior parameters $\sigma_{l}$ satisfy, for some $\alpha
>1/2$ and $0<r\le\alpha-\frac{1}4$,
%
%e14 #&#
\begin{eqnarray}\label{cond-sil}
\sigma_{l}&\ge& 2^{-l(\alpha+1/2)}\qquad\mbox{(log-Lipschitz case)},
\nonumber\\[-8pt]\\[-8pt]\nonumber
\sigma_{l}&=& 2^{-l(1/2+r)}\qquad\mbox{(Gaussian-case)}.
\end{eqnarray}
Typically, see examples below, such priors $f$ in (\ref
{def-prior-log}) under $\varphi=\varphi_G$ or $\varphi_H$ and (\ref
{cond-sil}) attain the rate $\bar\varepsilon_{n,\alpha}$ in (\ref
{strate}) in terms of Hellinger loss, up to logarithmic\vadjust{\goodbreak} terms. For some
$\nu>0$, suppose
%
%e15 #&#
\begin{equation}
\label{rate-hell} E_{f_0}^n\Pi\bigl[f\dvtx  h(f,f_0)>
(\log n)^\nu\bar\varepsilon_{n,\alpha} \given X^{(n)}
\bigr] \to0.
\end{equation}
If (\ref{rate-hell}) holds for some $\nu>0$, we denote $\varepsilon
_n:=(\log n)^\nu\bar\varepsilon_{n,\alpha}$
and $\zeta_n:=\varepsilon_n 2^{L_n/2}$, with $L_n$ as in (\ref{def-jn}).
%
%th2 #&#
\begin{teo} \label{teo-logdens}
Consider observations $X^{(n)}$ from model (\ref{mod-density}).
Suppose $\log f_0$ belongs to $ {\mathcal{C}}^{\alpha}[0,1]$, with
$\alpha\ge1$.
Let $\Pi$ be the prior on $ {\mathcal{F}}$ defined by (\ref
{def-prior-log}), with $\varphi=\varphi_G$ or $\varphi_H$. Suppose
that $\sigma_{l}$ satisfy (\ref{cond-sil}) and that (\ref
{rate-hell}) holds. Then, for $\alpha>1$ and $\varepsilon_{n,\alpha
}^*$ defined by (\ref{mmrate}), any $M_n\to\infty$, it holds, as
$n\to\infty$,
\[
E_{f_0}^n \Pi\bigl[f\dvtx  \llVert f-f_0\rrVert
_{\infty} > M_n \varepsilon_{n,\alpha}^* \given
X^{(n)} \bigr] \to0.
\]
In the case $\alpha=1$, the same holds with $\varepsilon_{n,\alpha
}^*$ replaced by $(\log{n})^{\eta}\varepsilon_{n,\alpha}^*$, for
some $\eta>0$.
\end{teo}
%
%The previous result holds under the more general condition $\log f_0
%slightly stronger H\"{o}lder condition is needed in the case $\al=1$
%only.

Theorem~\ref{teo-logdens} implies that log-density priors for many
natural priors on the coefficients achieve the precise optimal minimax
rate of estimation over H\"{o}lder spaces under sup-norm loss, as soon
as the regularity is at least $1$.

In\vspace*{1pt} the case $1/2<\alpha<1$, examination of the proof reveals that the
presented techniques %immediately
provide the sup-norm rate $\rho_n=n^{(1/2-3\alpha/2)/(1+2\alpha)}$ %$\rho_n = n^{\frac{1-2\al}{1+2\al}}$
up to logarithmic terms.
%One can note that, still up to log-factors, it holds $\rho_n=
For $1/2<\alpha<1$, we have $\varepsilon_{n,\alpha}^*\ll\rho_n \ll
\zeta_n$. So, although the minimax rate is not exactly attained for
those low regularities, %it may still be worth mentioning that
the obtained rate improves on the intermediate rate $\zeta_n$, which
was obtained in \cite{gn11} for slightly different priors. In the next
subsection, a prior is proposed which attains the minimax rate for the
sup-norm in the case $1/2<\alpha<1$.

Let us give some examples of prior distributions satisfying the
assumptions of Theorem~\ref{teo-logdens}.
In the Gaussian case, any sequence of the type $\sigma_{l}=2^{-l(1/2+\gamma)}$ with $0< \gamma\le\alpha-1/4$ satisfies both (\ref
{cond-sil}) and (\ref{rate-hell}). In the log-Lipschitz case, the
choice $\varphi=\varphi_{H,\tau}$ in (\ref{phi-ht}) with any $0\le
\tau<1$ combined with $\sigma_{l}=2^{-l\alpha}$ satisfies (\ref
{cond-sil})--(\ref{rate-hell}). Both claims follow from minor
adaptations of Theorem 4.5 in \cite{vvvz} and Theorem~2.1 in \cite
{rr12}, respectively; see Lemma~\ref{lem-exa}. In both Gaussian and
log-Lipschitz cases, we in fact expect (\ref{rate-hell}) to hold true
for many other choices of $\sigma_{l}$ under (\ref{cond-sil}) and
$\log\varphi_H$ Lipschitz, or under $\sigma_l\ge2^{-l(1/4+\alpha
)}$ in the Gaussian case, although such a general statement
in Hellinger distance is not yet available in the literature, to the
best of our knowledge.

%s3.2.2 #&#
\subsubsection{Random dyadic histograms}

Associated to the regular dyadic partition of $[0,1]$ at level $L \in
\mathbb{N}^*$, given by
$I_0^L=[0,2^{-L}]$ and $I_k^L=(k2^{-L},(k+1)2^{-L}]$ for $k=1,\ldots,2^L-1$, is a natural notion of histogram
\[
{\mathcal{H}}_L =\Biggl\{h\in L^\infty[0,1], h(x) =
\sum_{k=0}^{2^L-1} h_k \mathbh{1}_{I_k^L}(x),
h_k \in\mathbb{R}, k=0,\ldots,2^L-1 \Biggr\}\vadjust{\goodbreak}
\]
the set of all histograms with $2^L$ regular bins on $[0,1]$. Let
$\mathcal S_L = \{ \omega\in[0,1]^{2^L}; \sum_{k=0}^{2^L-1} \omega
_k =1\}$
be the unit simplex in $\mathbb{R}^{2^L}$. Further denote
\[
{\mathcal{H}}_L^1=\Biggl\{f\in L^\infty[0,1],
  f(x) = 2^L \sum_{k=0}^{2^L-1}
\omega_k \mathbh{1}_{I_k^L}(x), (\omega
_0,\ldots,\omega_{2^L-1})\in{\mathcal{S}}_L
\Biggr\}.
\]
The set $ {\mathcal{H}}_L^1$ is the subset of $ {\mathcal
{H}}_L$ consisting of histograms which are densities on~$[0,1]$. Let
$ {\mathcal{H}}^1$ be the set of all histograms which are densities
on $[0,1]$.
%The set $\cH_L$ is a closed subspace of $L^2[0,1]$. For any function
%$h$ in $L^2[0,1]$, consider its projection $h_{[L]}$ in the

A simple way to specify a prior on $ {\mathcal{H}}_L^1$ is to set
$L=L_n$ deterministic and to fix a distribution
for $\omega_L:=(\omega_0,\ldots,\omega_{2^L-1})$.
Set $L=L_n$ as defined in (\ref{def-jn}).
Choose some fixed constants $a, c_1,c_2>0$ and let
%and set, for $L_n=\lfloor\log_2 J_n \rfloor$ and $J_n= \lfloor(n/
%
%e16 #&#
\begin{equation}
\label{prk} L = L_n,\qquad\omega_L \sim{\mathcal{D}}(
\alpha_0, \ldots, \alpha_{2^L-1}), \qquad c_1
2^{-La}\le\alpha_k \le c_2,
\end{equation}
for any admissible index $k$, where $ {\mathcal{D}}$ denotes the
Dirichlet distribution on $S_L$. Unlike suggested by the notation, the
coefficients $\alpha$ of the Dirichlet distribution are allowed to
depend on $L_n$, so that $\alpha_{k}=\alpha_{k,L_n}$. % For
%simplicity of notation, we keep the $\al_k$ notation.
%Given $L=L_n$, denote by $\cV_{L_n}$ the space spanned by the elements
%of the Haar basis up to level $L_n-1$, that is by $\varphi^H,
%Given $L=L_n$, denote by $H_{L_n}$ the collection of the Haar basis
%elements up to level $L_n-1$, that is by $H_{L_n}=\{ \varphi^H,
%by $\ga_k$ the (constant) value of $\ga$ on the interval $I_{k}^L$.
%We impose the following mild restriction on the $\alpha$'s.
%For any $\ga\in H_{L_n}$,
%Suppose
%o\Big( \sqrt{{n}\log{n}} \Big).
% petit o ou grand o ?
%
%th3 #&#
\begin{teo} \label{teo-histo}
Let $f_0\in{\mathcal{F}}_0$ and suppose $f_0$ belongs to $
{\mathcal{C}}^{\alpha}[0,1]$, where $1/2 < \alpha\le1$.
Let $\Pi$ be the prior on $ {\mathcal{H}}^1\subset{\mathcal
{F}}$ defined by (\ref{prk}). %Suppose \eqref{cond-dir} holds.
Then, for $\varepsilon_{n,\alpha}^*$ defined by (\ref{mmrate}) and
any $M_n\to\infty$ it holds,
as $n\to\infty$,
\[
E_{f_0}^n \Pi\bigl[f\dvtx  \llVert f-f_0\rrVert
_{\infty} > M_n \varepsilon_{n,\alpha}^* \given
X^{(n)} \bigr] \to0.
\]
\end{teo}
According to Theorem~\ref{teo-histo}, random dyadic histograms achieve
the precise minimax rate in sup-norm over H\"{o}lder balls.
Condition~(\ref{prk}) is quite mild. For instance, the uniform choice
$\alpha
_0=\cdots=\alpha_{2^L-1}=1$ is allowed, as well as a variety of
others, for instance, one can take $\alpha_k=\alpha_{k,L_n}$ to
originate from a measure $A=A_{L_n}$ on the interval $[0,1]$, of finite
total mass $\bar A_{L_n}:=A([0,1])$. By this we mean, $\alpha
_k=A(I_k^{L_n})$. If $A/\bar A_{L_n}$ has say a fixed continuous and
positive density $a$ with respect to Lebesgue measure on $[0,1]$, then
(\ref{prk}) is satisfied as soon as there exists a $\delta>0$ with
$2^{-\delta L_n}\lesssim\bar A_{L_n}\lesssim2^{L_n}$.

%Given \eqref{prk}, it essentially asks for the collection of $\al$'s
%to vary in a not too abrupt way along the indexes $k$. To gain
%insight, suppose that $\al_k=\al_{k,L_n}$ originate from a measure
%$A=A_{L_n}$ on the interval $[0,1]$, of finite total mass $\bar
%A_{L_n}:=A([0,1])$. By this we mean $\al_k=A(I_k^L)$. If $A/\bar
%A_{L_n}$ has say a fixed density $a$ with respect to Lebesgue measure,
%then \eqref{cond-dir}
% becomes
% \begin{equation} \label{cond-mes}
% \bar A_{L_n} \max_{\ga\in H_{L_n}}
% \Big\left|\int_0^1 \ga(u) a(u) \,du\Big\right| = o\Big( \sqrt{n\log{n}}
% \end{equation}
%For instance, \eqref{cond-dir} is satisfied for the classical choice $

%s3.2.3 #&#
\subsubsection{Further examples}
A referee of the paper, whom we thank for the suggestion, has asked
whether the proposed technique would work for other priors, more
specifically for non-$n$-dependent priors in density estimation.
Although not considered here for lack of space, we would like to
mention the important class of P\'olya tree priors; see, for example,
Lavine \cite{lavine92}.
For well-chosen parameters, it can be shown that these priors achieve
supremum-norm consistency in density estimation (consistency in the,
weaker, Hellinger sense was studied, e.g., in \cite{bsw99}) and
minimax rates of convergence in the sup-norm can be obtained. In
particular, this class contain canonical (i.e., non-$n$-dependent)
priors that achieve such optimal rates in density estimation. This will
be studied elsewhere.

%s3.3 #&#
\subsection{Discussion}

We have introduced new tools which allow to obtain optimal minimax
rates of contraction in strong distances for posterior distributions.
The essence of the technique is to view the problem semiparametrically
as the uniform study of a collection of semiparametric Bayes
concentration results, very much in the spirit of nonparametric
Bernstein--von Mises results as studied in \cite{CN12}.
For the sake of clarity, we refrain of carrying out further extensions
in the present paper but briefly mention a few applications.
%one can note that the impact of the techniques goes fairly much beyond
%this. Let us mention a few examples.
From the sup-norm rates, optimal results---up to logarithmic terms- in
$L^q$-metrics, $q\ge2$, can be immediately obtained by interpolation.
Adaptation to the unknown $\alpha$ could also be considered. This will
be the object of future work. However, note that ``fixed $\alpha$''
nonparametric results as such are already very desirable in strong
norms. They can, for instance, be used in the study of remainder terms of
semiparametric functional expansions or of LAN-expansions as, for
example, to check the conditions of application of semiparametric
Bernstein--von Mises theorems as in \cite{ic12}. In this semiparametric
perspective,
%Although interesting from the nonparametric perspective,
adaptation to $f$ is in fact not always desirable, %the semiparametric
%point of view of studying functionals
since posteriors for functionals may behave pathologically when an
adaptive prior on the nuisance is chosen; see \cite{rr12} and \cite
{cr13}, where it is shown that too large discrepancies in smoothness
between the semiparametric functional and the unknown $f$ can lead to
undesirable bias. %Nevertheless, obtaining results for adaptive priors
%with the present method is conceivable and is the subject of current
%research.
Also, we expect the present methodology to give results in a broad
variety of statistical models and/or for different classes of priors.
Indeed, it reduces the problem of the strong-distance rate to two
parts: (1)~uniform semiparametric study of functionals and (2)~high-frequency
bias. The first part is very much related to obtaining (uniform)
semiparametric Bernstein--von Mises (BvM) results. So, any advance in
BvM theory for classes of priors will automatically lead to advances
in~(1). As for~(2),
the studied examples suggest that for frequencies above
the cut-off the posterior behaves essentially as the prior itself. So,
contrary to the BvM-regime (1) where the prior washes out
asymptotically, one does not expect a universal behaviour for this
part. However, showing that the posterior is close to the prior
provides a possible method of proof.
%although ad-hoc methods can be used depending on the model and prior
%at hand (e.g. choosing priors cutting high-frequencies or applying a
%concentration of measure inequality), the example of log-density
%priors indicates that a pre-processing step may be useful in general
%to `adapt' the problem first to the `geometry' of the prior.%: for
%log-densities, the logarithmic transformation exactly plays this role.

%adaptive histogram ?) or nonparametric BvM's in density estimation ?
%
% Mention that we could presumably use localised needlet frames to
%transport results to other geometries -at least already in the
%nonadaptive case-?
%
% Lower bound for log-density series prior (say in the `f_0 more
%regular' case ?)

%s4 #&#
\section{Proofs}

%s4.1 #&#
\subsection{Gaussian white noise}

%le1 #&#
\begin{lem} \label{lem-ratio}
Let $X^{(n)}$ follow model (\ref{gwn}). Let $f_0$ satisfy (\ref
{tolk}) and let the prior $\Pi$ be chosen according to (\ref{prior-unif}).
There exists $C>0$ such that for any real $t$, any $n\ge2$ and $l\le
L_n$, with $L_n$ defined in (\ref{def-jn}),
\[
E_{f_0}^nE^{\Pi}\bigl[ e^{t \sqrt{n}(f_{lk} - x_{lk})} \given
X^{(n)}\bigr] \le C e^{t^2/2}.
\]
\end{lem}
\begin{pf}% {Proof of Lemma~\ref{lem-ratio}}
The proof is similar to the first lines of the proof of Theorem~5 in
\cite{CN12}: one uses Bayes' formula to express the posterior
expectation in the lemma. Next, using (\ref{def-jn}) and\vadjust{\goodbreak} (\ref
{tolk}), one checks that for any $v\in[-L_0,L_0]$ with $L_0:= \sqrt
{2}(B-R)/2$, the ratio $\llvert f_{0,lk}+v/\sqrt{n}\rrvert/\sigma_{l}$
is at most
$R+(B-R)/2<B$, for any $l\le L_n$ and $k$.
For such $v$'s, since $\varphi$ is the uniform density on $[-B,B]$,
the expression involving $\varphi$ in the next line is constant,
%%equal to $1/(2B)$,
and thus can be removed from the expression, leading~to
% finds a universal $L_0>0$
%such that for any $l,k$ and $v\in(-L_0,L_0)$, the expression involving
%the density $\varphi$ in the next line is bounded below (and also
%above}) by a fixed constant, and thus can be removed from the
%expression, leading to
%
\begin{eqnarray*}
&& E^{\Pi}\bigl[ e^{t \sqrt{n}(f_{lk} - x_{lk})} \given X^{(n)}\bigr]
\\
&&\qquad  =
e^{-t\varepsilon_{lk}} \frac{ \int e^{-v^2/2+(t+\varepsilon_{lk})v}
\varphi((f_{0,lk} + v/\sqrt{n})/\sigma_l ) \,dv}{
\int e^{-v^2/2+\varepsilon_{lk}v} \varphi((f_{0,lk} +
v/\sqrt{n})/\sigma_l) \,dv }
\\
&&\qquad  \lesssim e^{-t\varepsilon_{lk}} \frac{ \int e^{tv-(v-\varepsilon
_{lk})^2/2} \,dv}{\int_{-L_0}^{L_0} e^{ -(v-\varepsilon_{lk})^2/2} \,dv}
\lesssim\frac{ \int e^{tu-u^2/2} \,du}{\int_{-L_0}^{L_0} e^{
-(v-\varepsilon_{lk})^2/2} \,dv}
\\
&&\qquad  \lesssim e^{t^2/2} \biggl[ \int_{-L_0}^{L_0}
e^{ -(v-\varepsilon
_{lk})^2/2} \,dv \biggr]^{-1}.
\end{eqnarray*}
Since $\varepsilon_{lk}$ are standard normal, simple calculations show
that the expectation of the inverse of the quantity under brackets is
bounded by a universal constant, as in \cite{CN12}, pages 2015--2016.
%, which concludes the proof.
%Therefore, for any real $t$ and $l\le J_n$, as in \cite{CN12}, which
%leads to
\end{pf}

\begin{pf*}{Proof of Theorem~\ref{teo-gwn}}
\emph{Small $l$}. Let\vspace*{1pt} us first consider indexes $l$ with $l\le L_n$. For
any real $t$, set $Q_{lk}(t):= E_{f_0}^n E^{\Pi}[ e^{t\sqrt
{n}(f_{lk}-X_{lk})} \given X^{(n)} ]$.
Using the fact that $\varphi$ is bounded,
\[
Q_{lk}(t) \lesssim E_{f_0}^n
\frac{\int e^{t(v-\varepsilon _{lk})-((v-\varepsilon_{lk})^2/2)} \,dv}{ \int e^{-(v-\varepsilon_{lk})^2/2} \varphi((f_{0,lk}+v/\sqrt{n})/\sigma_{l} )\,dv }.
\]
Introduce the set, for any possibly $l$-dependent sequence $M_l$,
%
%e17 #&#
\begin{equation}
\label{seta} {\mathcal{A}}(M_l):= \biggl\{ v\dvtx  \biggl\llvert
\frac{f_{0,lk}+v/\sqrt{n}}{\sigma_{l}} \biggr\rrvert\le M_l \biggr\}.
\end{equation}
Choose $M_l=C (l+1)^{\mu}$ with $\mu=(1+\delta)^{-1}$. This implies,
with our choices of $M_l, \sigma_l$ and taking $C$ large enough, that
$ {\mathcal{A}}(M_l)$ contains the interval $(-1,1)$.
%First one restricts the integral on the denominator to $(-1,1)$ and
%next notes that for any $v\in(-1,1)$, any $l\le L_n$
First restricting the integral on the denominator to $(-1,1)$ and next
using the tail condition on $\varphi$ and the fact that $\varphi\ge
c_\varphi>0$ on $(-1,1)$, one gets
\begin{eqnarray*}
Q_{lk}(t) & \lesssim& E_{f_0}^n
\frac{e^{t^2/2} }{ e^{-l}
\int_{-1}^1 e^{-(v-\varepsilon_{lk})^2/2} \,dv } \lesssim e^{t^2/2+l}.
\end{eqnarray*}
The maximal inequality argument from Section~\ref{firstex}
directly yields $(\mathrm{i})\le\varepsilon_{n,\alpha}^*$.

\emph{Large $l$}. Let us now consider the case $l>L_n$. For any real $t$ set,
\begin{eqnarray*}
&& E_{f_0}^n E^{\Pi}\bigl[ e^{t f_{lk}} \given
X^{(n)} \bigr]
\\
&&\qquad = E_{f_0}^n \frac{ \int e^{t(f_{0,lk}+v/{\sqrt
{n}})} e^{-v^2/2+\varepsilon_{lk}v} 1/{\sqrt{n}\sigma_{l}}
\varphi( (f_{0,lk}+v/{\sqrt{n}})/{\sigma_{l}}
)\,dv }{ \int e^{-v^2/2+\varepsilon_{lk}v} 1/{\sqrt{n}\sigma_{l}}
\varphi( ({f_{0,lk}+v/{\sqrt{n}}})/{\sigma_{l}} ) \,dv }
\\
&&\qquad  =:
E_{f_0}^n \frac{N_{lk}(t)}{D_{lk}}.
\end{eqnarray*}
To bound the denominator, first restrict the integral to the set $
{\mathcal{A}}:= {\mathcal{A}}(1)$ as defined in (\ref{seta}). Set
\[
\zeta_l=\int_{ {\mathcal{A}}} v \frac{1}{\sqrt{n}\sigma_{l}} \varphi
\biggl( \frac{f_{0,lk}+v/{\sqrt{n}}}{\sigma_{l}} \biggr)\,dv,
\]
next apply Jensen's inequality with the logarithm function to get, with
$\llvert{\mathcal{A}}\rrvert$ the diameter of $ {\mathcal{A}}$ and some
constant $C>0$,
\begin{eqnarray*}
\log D_{lk} & \ge&-\frac{ \llvert{\mathcal{A}}\rrvert\llVert\varphi
\rrVert_\infty
}{2\sqrt{n}\sigma_{l}} \sup_{v\in{\mathcal{A}}}v^2
+ \varepsilon_{lk}\zeta_l
\\
& \ge&- Cn\bigl(\sigma_{l}^2 + f_{0,lk}^2
\bigr) + \varepsilon_{lk}\zeta_l,%
\end{eqnarray*}
where we have used that $M_l=1$ in (\ref{seta}). Below we shall also
use that
%The following bound will be useful in a moment
%
\begin{eqnarray*}
\llvert\zeta_l\rrvert& \lesssim&\frac{ \llvert{\mathcal{A}}\rrvert
\llVert\varphi\rrVert_\infty
}{\sqrt{n}\sigma_{l}} \sup
_{v\in{\mathcal{A}}}\llvert v\rrvert\lesssim\sqrt{n}\bigl(\llvert
f_{0,lk}\rrvert+\sigma_{l}\bigr).% \leqa\rn\sil.
\end{eqnarray*}
To\vspace*{1pt} bound the numerator from above, split the integrating set
into $ {\mathcal{A}}:= {\mathcal{A}}(1)$ and~$ {\mathcal
{A}}^c$ and
write $N_{lk}(t)=: N_{lk}^{(1)}(t) + N_{lk}^{(2)}(t)$ for the integrals
over each respective set.
Using the previous bound on $D_{lk}$, that $\llvert t(f_{0,lk}+\frac
{v}{\sqrt
{n}})\rrvert\le\llvert t\rrvert\sigma_{l}$ by definition of~$
{\mathcal{A}}$, and
Fubini's theorem,
\begin{eqnarray*}
&& E_{f_0}^n \frac{N_{lk}^{(1)}(t)}{D_{lk}}
\\
&&\qquad  \le e^{\llvert t\rrvert
\sigma_{l}+ Cn(\sigma_{l}^2 + f_{0,lk}^2)} \int
_{ {\mathcal{A}}} E_{f_0}^n \bigl[ e^{(v-\zeta_l)\varepsilon_{lk}}
\bigr] \frac{e^{-v^2/2}}{\sqrt{n}\sigma_{l}} \varphi\biggl(
\frac
{f_{0,lk}+v/{\sqrt{n}}}{\sigma_{l}} \biggr)\,dv
\\
&&\qquad  \le e^{\llvert t\rrvert\sigma_{l}+ Cn(\sigma_{l}^2 + f_{0,lk}^2) +
\zeta _l^2/2}\llVert\varphi\rrVert_{\infty} \sup
_{v\in{\mathcal{A}}} e^{\llvert v\zeta_l\rrvert}
\\
&&\qquad  \lesssim e^{\llvert t\rrvert\sigma_{l}+ C'n(\sigma_{l}^2 +
f_{0,lk}^2)}. %\leqa
%e^{t\sil+ Cn(\sil^2 f_{0,lk}^2)}
\end{eqnarray*}
On the other hand, the term over $ {\mathcal{A}}^c$ can be bounded
as follows:
\begin{eqnarray*}
&& E_{f_0}^n \frac{N_{lk}^{(2)}(t)}{D_{lk}}
\\
&&\qquad  \lesssim  e^{ Cn(\sigma
_{l}^2 + f_{0,lk}^2) }
\\
&&\quad\qquad{}\times
\int_{(-1,1)^c} e^{t\sigma_{l}w} E_{f_0}^n
\bigl[ e^{-n/2 (w\sigma_{l}-f_{0,lk})^2+\varepsilon_{lk}(\sqrt{n}(w\sigma
_{l}-f_{0,lk})-\zeta_l)} \bigr] \varphi(w)\,dw
\\
&&\qquad  \lesssim e^{ Cn(\sigma_{l}^2 + f_{0,lk}^2) + \zeta_l^2/2
}\int_{(-1,1)^c} e^{t\sigma_{l}w-\sqrt{n}\zeta_l (\sigma_{l}w - f_{0,lk})}
\varphi(w)\,dw
\\
&&\qquad  \lesssim e^{ C'n(\sigma_{l}^2 + f_{0,lk}^2) } \int_{(-1,1)^c}
e^{(t\sigma_{l}+\sqrt{n}\sigma_{l}\zeta_l) w}
\varphi(w)\,dw.
\end{eqnarray*}
Using the tail behaviour of $\varphi$ leads to
\begin{eqnarray*}
E_{f_0}^n \frac{N_{lk}^{(2)}(t)}{D_{lk}} & \lesssim&
e^{ C'n(\sigma_{l}^2 + f_{0,lk}^2) }e^{C \{ \sigma_l(\llvert t\rrvert
+\sqrt
{n}\llvert\zeta_l\rrvert)\}^{(\delta+1)/{\delta}}}.
\end{eqnarray*}
One deduces, using that for $l>L_n$, one has $n(\sigma_{l}^2 +
f_{0,lk}^2)\le n 2^{-l(1 +2\alpha)} \lesssim \log{n}\lesssim l$, that
for $t>0$,
\begin{eqnarray*}
R_{lk}(t)&:=& E_{f_0}^nE^\Pi\Bigl[
\max_{k}\llvert f_{lk}\rrvert\big\given
X^{(n)}\Bigr]
\\
 & \lesssim&\frac{1} t \bigl(l + \log
\bigl(e^{t\sigma_l}+e^{C \{ \sigma_l(t+\sqrt
{n}\llvert\zeta_l\rrvert)\}^{(\delta+1)/\delta} }\bigr)\bigr)
\\
& \lesssim&\frac{l}{t} + \sigma_{l}+ \frac{1}{t} \bigl\{
\sigma_l\bigl(t+\sqrt{n}\llvert\zeta_l\rrvert\bigr)
\bigr\}^{(\delta+1)/\delta}.
\end{eqnarray*}
Set $t=\sigma_l^{-1}l^{\delta/(\delta+1)}$ to deduce, using
$\sigma_{l}\sqrt{n}\llvert\zeta_l\rrvert\lesssim l^{\delta/(\delta+1)}$
for $l>L_n$,
\[
R_{lk}(t) \lesssim l^{1/(\delta+1)}\sigma_l \lesssim
2^{-l(1/2 + \alpha)}
\]
and further obtain $\mathrm{(ii)} \le\sum_{l>L_n} 2^{l/2} 2^{-l(1/2 + \alpha)} \lesssim h_n^{\alpha}
=\varepsilon_{n,\alpha}^*$. Therefore, for any \mbox{$\delta>0$}, the rate
is precisely $\varepsilon_{n,\alpha}^*$.
\end{pf*}

%s4.2 #&#
\subsection{Density estimation, notation} \label{sec-dens-not}
%%roadmap of the proofs}

%The purpose of this section is first to introduce notation useful in
%the proofs for density estimation.
%Next, we present the main argument in the simplest case of random
%dyadic histogram, where the approximation part is simplified by the
%favorable properties of approximation by histograms
%in combination with the Haar basis.

%{\em Notation density case.}
Given observations $X^{(n)}$ from (\ref{mod-density}), denote by $\ell
_n(f)$ the log-likelihood $ \ell_n(f) = \sum_{i=1}^n \log f(X_i)$. %,
%well-defined for any density $f$ in the class $\cF$.
For any $u,v$ in $L^2(P_{f_0})=:L^2(f_0)$, define the inner-product
${\langle}\cdot,\cdot{\rangle}_L$ with associated norm $\llVert\cdot
\rrVert
_L$, together with a stochastic term
$W_n(u)$, as follows:
\begin{eqnarray*}
{\langle}u, v {\rangle}_L & =& \int_0^1
(u - P_{f_0}u) (v-P_{f_0}v) f_0,
\\
W_n(u) & =& \frac{1}{\sqrt{n}}\sum
_{i=1}^n \bigl[u(X_i)-P_{f_0}u
\bigr].
\end{eqnarray*}
In particular, in empirical process notation $W_n(u)=\mathbb{G}_n(u)$.
For any $f$ in $ {\mathcal{F}}$, set $R_n(f,f_0)=\sqrt
{n}P_{f_0}\log(f/f_0) + n\llVert\log(f/f_0) \rrVert_L^2/2$.
For any $f\in{\mathcal{F}}$, it holds
\[
\ell_n(f) - \ell_n(f_0) = -
\frac{n}{2}\bigl\llVert\log(f/f_0) \bigr\rrVert
_L^2 + \sqrt{n}W_n\bigl( \log(f /
f_0) \bigr) + R_n(f,f_0).
\]
Denote, for any density $f$ in $ {\mathcal{F}}$ and any given $u$
in $L^2(f_0)$,
\begin{eqnarray*}
\mathcal{B}(u,f,f_0) & =& {\biggl\langle}\frac{f-f_0}{f_0}, u {
\biggr\rangle}_L - {\bigl\langle}\log(f/f_0), u {\bigr
\rangle}_L.
\end{eqnarray*}
%
%R_n(\eta,\eta_0) = \ell_n(\eta)-\ell_n(\eta_0)
%+\frac{1}{2} \left\Vert h \right\Vert_L^2 - W_n(h) +,
Let $D_n$ be a measurable set. Denote by $\Pi^{D_n}$ the restriction
of $\Pi$ to $D_n$. Suppose, as $n\to\infty$,
%
%e18 #&#
\begin{equation}
\label{set} E_{f_0}^n \Pi\bigl(D_n \given
X^{(n)}\bigr) =1+o(1).
\end{equation}
%
%For any $M>0$, Markov's inequality implies
%X^{(n)}] \le(M
%E_{f_0}^n \int_{D_n} \left\Vert f-f_0\right\Vert_\infty d\Pi(f\given
%X) + o(1). \]
%Denote by $\Pi^{D_n}$ the prior distribution conditioned to lie in the
%set $D_n$.
Combining (\ref{set}) and Markov's inequality leads to, for any
$M_n\to\infty$,
\begin{eqnarray*}
&& E_{f_0}^n \Pi\bigl[f\dvtx  \llVert f-f_0
\rrVert_\infty>M_n \varepsilon_{n,\alpha}^* \given
X^{(n)}\bigr]
\\
&&\qquad  \le\bigl(M_n\varepsilon_{n,\alpha}^*\bigr)^{-1}
E_{f_0}^n \bigl[ E^{\Pi^{D_n}}\bigl[ \llVert
f-f_0\rrVert_\infty\given X^{(n)}\bigr]\Pi
\bigl(D_n \given X^{(n)}\bigr) \bigr]+ o(1).
\end{eqnarray*}
In the sequel, we focus on bounding $E^{\Pi^{D_n}}[ \llVert
f-f_0\rrVert_\infty
\given X^{(n)} ]$ from above.
%shows that for the posterior to concentrate in a
%$\left\Vert\cdot\right\Vert_\infty$-ball of radius $M_n\veps_{n,
% is enough to prove that there exists a finite $C>0$ such that
%E_{f_0}^n \int_{D_n} \left\Vert f-f_0\right\Vert_\infty d\Pi(f\given
%X) \le C\veps_{n,

%s4.3 #&#
\subsection{Density estimation, log-density priors} \label{sec-dens-proof}

Let us define the set $D_n$ by, for $\varepsilon_n=(\log n)^\nu\bar
\varepsilon_{n,\alpha}$ the rate in
(\ref{rate-hell}), $L_n$ as in (\ref{def-jn}) and $\zeta
_n=\varepsilon_n 2^{L_n/2}$,
%
%e19 #&#
\begin{equation}
\label{def-dn} D_n = \bigl\{f, \llVert f-f_0\rrVert
_2\le\varepsilon_n, \llVert f-f_0\rrVert
_{\infty}\le\zeta_n \bigr\}.
\end{equation}
It follows from Lemma~\ref{lem-inter} below that $\Pi(D_n\given
X^{(n)})$ goes to $1$ in probability, up to replacing $\varepsilon_n$
by $M\varepsilon_n$ for a large enough constant $M$, and similarly for
$\zeta_n$. Indeed, since a $\varepsilon_n$-Hellinger-contraction rate
for the posterior is assumed, see (\ref{rate-hell}), the conditions of
Lemma~\ref{lem-inter} are satisfied.
%Under \eqref{rate-hell}, Lemma~\ref{lem-inter} implies \eqref{set}, up
%to multiplication of $\veps_n,\zeta_n$ by large enough constants.

%s4.3.1 #&#
\subsubsection{First step, reduction to the logarithmic scale}

Let us set $g=\log f$ and $g_0=\log f_0$. With $T$ defined in (\ref
{def-T}), one has $g=T-c(T)$. First, one notes that obtaining a rate
going to $0$ for $\llVert g-g_0\rrVert_\infty$ implies the same rate
up to
constants for $\llVert f-f_0\rrVert_\infty$. Indeed, $\llVert
f-f_0\rrVert_\infty=\llVert
e^{g_0}(e^{g-g_0}-1)\rrVert_\infty
\lesssim\llVert g-g_0\rrVert_\infty$ using the bound $\llvert
e^x-1\rrvert\lesssim\llvert x\rrvert$ for
small $x$ and that $\llVert f_0\rrVert_\infty$ is bounded. So,
instead of writing
Markov's inequality as above with $f$, we write it with $g$, the set
$D_n$ still being the one defined in (\ref{def-dn}) with the
dependence on $f-f_0$.

That is, we focus on bounding $E^{\Pi^{D_n}}[ \llVert g-g_0\rrVert
_\infty\given
X^{(n)} ]$ from above.
Now write, with the notation $g^{L_n}$ denoting the $L^2$-projection up
to level $L_n$ as in Section~\ref{firstex}, and $L_n$ as in (\ref{def-jn}),
%The reason we work at exponential level is to keep things relatively
%simple to handle the {\it bias} (in the context of this paper this
%means the terms with large $j$).
%
\begin{eqnarray*}
&& E^{\Pi^{D_n}}\bigl[ \llVert g-g_0\rrVert_{\infty}
\given X^{(n)} \bigr]
\\
&&\qquad  \le\underbrace{ \int\bigl\llVert g^{L_n} - g_0^{L_n}
\bigr\rrVert_{\infty} \,d\Pi^{D_n}\bigl(f\given X^{(n)}
\bigr)}_{(\mathrm{i})} + \underbrace{\int\bigl\llVert g^{L_n^c}\bigr
\rrVert_{\infty} \,d\Pi^{D_n}\bigl(f\given X^{(n)}
\bigr)}_{(\mathrm{ii})} + \underbrace{\bigl\llVert g_0^{L_n^c}
\bigr\rrVert_{\infty}}_{(\mathrm{iii})}.
\end{eqnarray*}
The term (ii) is $0$ because the sum defining $T$ goes up to level $l
\le L_n$ under the prior distribution, and the constant function $1$ is
orthogonal to higher levels. Since $g_0=\log f_0$ belongs to $B_{\infty,\infty}^{\alpha}$ by assumption, the term (iii) is bounded
by a constant times $\varepsilon_{n,\alpha}^*$.

We now start analysing the term (i). First, let us introduce, for $\{
{\mathcal{A}}_{l,k}\}_{l,k}$ a collection of
elements of $L^2(f_0)$ to be chosen later, and $L_n$ as in (\ref{def-jn}),
%
%e20 #&#
\begin{equation}
\label{sto} \Gamma^{L_n}(\cdot):= g_0^{L_n}(
\cdot) + \frac{1}{\sqrt{n}} \sum_{l=0}^{L_n}
\sum_{k=0}^{2^l -1} W_n( {
\mathcal{A}}_{l,k} ) \psi_{lk}(\cdot).
\end{equation}
Next, let us write
\[
\mbox{(i)} \le\int\bigl\llVert g^{L_n} - \Gamma^{L_n}\bigr\rrVert
_{\infty} \,d\Pi^{D_n}\bigl(f\given X^{(n)}\bigr) + \bigl
\llVert\Gamma^{L_n}-g_0^{L_n}\bigr\rrVert
_{\infty}.
\]
The second term is bounded with the help of Lemma~\ref{lem-Gamma}.
For the first term, following the scheme of proof of the maximal
inequality in Section~\ref{firstex} via the moment generating
function, one sees that it is enough to bound for $t>0$ the following
quantity, uniformly in $l,k$ and $l\le L_n$
%
%e21 #&#
\begin{equation}
\label{emlk} {\mathcal{M}}_{lk}(t):= e^{ -t W_n( {\mathcal{A}}_{l,k}) }
E^{\Pi^{D_n}} \bigl[ e^{t\sqrt{n}{\langle}g-g_0,\psi_{lk} {\rangle
}_2} \given X^{(n)} \bigr].
%e^{-tW_n()}.
\end{equation}
Denote $\rho(x):=\log(1+x)-x$. It holds
\begin{eqnarray*}
\int_0^1 (g-g_0)
\psi_{lk}& =& \int_0^1 \log\biggl[
\frac
{f-f_0}{f_0} + 1 \biggr] \psi_{lk}
\\
& =& \int_0^1 \frac{f-f_0}{f_0}
\frac{\psi_{lk}}{f_0} f_0 + \int_0^1
\rho\biggl(\frac{f-f_0}{f_0} \biggr) \psi_{lk}.
\end{eqnarray*}
On $D_n$ we have an intermediate sup-norm rate $\zeta_n=o(1)$ when
$\alpha>1/2$. In this case the argument of $\rho$ in the previous
display tends to $0$. Using the bound \mbox{$\llvert\rho(u)\rrvert\le u^2$}
for small
$u$, one gets
%
%e22 #&#
\begin{equation}
\label{te-rem} \biggl\llvert\int_0^1 \rho
\biggl(\frac{f-f_0}{f_0} \biggr)\psi_{lk}\biggr\rrvert\le\llVert\psi
_{lk}\rrVert_\infty\int_0^1
\biggl(\frac{f-f_0}{f_0} \biggr)^2 \lesssim2^{l/2} \llVert
f-f_0\rrVert_2^2.
\end{equation}
This bound is a $O(1/\sqrt{n})$ on $D_n$ as soon as
$2^{L_n/2}\varepsilon_n^2 =O(1/\sqrt{n})$, which is satisfied if
$\alpha> 1$.
% (quite interestingly, one can show that the present reduction at the
%logarithmic level is {\it not} what limits the rate in the case $1/2<
%What precedes suggests to write
%where $\doteq$ means up to a remainder term -here it has just been
%controlled above- and
This implies that the inner-product ${\langle}g-g_0,\psi_{lk}
{\rangle}_2$ in (\ref{emlk}) can be
replaced by ${\langle}f-f_0,\zeta_{l,k} {\rangle}$, where
%
%e23 #&#
\begin{equation}
\label{def-zet} \zeta_{l,k} = \frac{\psi_{lk}}{f_0}.
\end{equation}
That is, we can reason as if one would be considering the
semiparametric problem of
estimating the linear functional of the density
$f \to{\langle}\zeta_{l,k},f {\rangle}_2$. The corresponding
efficient influence function is
$\tilde\zeta_{l,k}=\zeta_{l,k}- P_{f_0}\zeta_{l,k}$, with respect
to the tangent set $ {\mathcal{H}}_{f_0}:=\{h\dvtx [0,1]\to\mathbb
{R}, h\mbox{ bounded, }\int_0^1 h f_0 = 0 \}$; see \cite{aad98},
Chapter~25 for definitions.

There is one difficulty with $\zeta_{l,k}$. It is not an element of
the basis of expansion of the prior
$\Pi$, %which in general prevents to use it directly for a change of
%variables: indeed, $\tilde\zeta_{l,k}$
so it needs to be properly approximated by the prior in some sense. In
fact, there is a fundamental difference with what has been done so far
in proving BvM-type results; see, for example, \cite{ic12,rr12,cr13}.
Here, we need to study approximating sequences \emph{uniformly} in the
indexes $l,k$ and a sharp control on this dependence is essential; see
the key Lemma~\ref{lem-approx}, where two regimes of indexes ``$l$''
arise, depending on whether $l$ is small or close to~$L_n$.

So, instead of working with $\zeta_{l,k}$ directly, one replaces it by
an approximation $ {\mathcal{A}}_{l,k}$ defined in (\ref{approx1})
below. %, leading to the scheme%, for a suitable $\zeta_n$
% \rn\psg f-f_0,\cA_{l,k} \psd_2 \]
This induces a \emph{bias} term for any $l,k$, familiar in the context
of semiparametric BvM results; see, for example, \cite
{ic12,ic122,cr13}, equal to
%
%e24 #&#
\begin{equation}
\label{sp-bias} \sqrt{n} {\langle}f-f_0,\zeta_{l,k}- {
\mathcal{A}}_{l,k} {\rangle}_2 = \sqrt{n}\int
_0^1 (f-f_0) (
\zeta_{l,k}- {\mathcal{A}}_{l,k} ).
\end{equation}
%
%and bounded above by $\left\Vert f-f_0\right\Vert_1\left\Vert
This term is controlled using Lemma~\ref{lem-approx} below. Indeed, on
$D_n$ the bounds of (\ref{sp-bias}) of Lemma~\ref{lem-approx} are at
most $\sqrt{n h_n}\varepsilon_n = o(1)$ if $\alpha>1$.
Next,
apply Lemma~\ref{lem-lapu} with $\gamma_n= {\mathcal{A}}_{l,k}$.
The estimates of $L^2$ and sup-norm of $ {\mathcal{A}}_{l,k}$ imply
that the conditions of application of Lemma~\ref{lem-lapu} are
satisfied. Thus,%One gets
%
%e25 #&#
\begin{equation}
\label{te-lap1} {\mathcal{M}}_{lk}(t) \le e^{Ct^2}
\frac{\int e^{ \ell_n(f_t) - \ell_n(f_0) } \,d\Pi^{D_n}(f) }{
\int e^{ \ell_n(f) - \ell_n(f_0) } \,d\Pi^{D_n}(f) },
\end{equation}
where we have set $f_t=e^{g_t}$ with $g_t=g_{t,l,k}$ defined as in
Lemma~\ref{lem-lapu} by (the expression is invariant under adding a
constant to $g$, so one can write it either with $ {\mathcal
{A}}_{l,k}$ or
$\tilde{\mathcal{A}}_{l,k}= {\mathcal{A}}_{l,k}-P_{f_0}
{\mathcal{A}}_{l,k}$)
%
%e26 #&#
\begin{equation}
\label{te-gt1} g_t = g - \frac{t}{\sqrt{n}} {\mathcal{A}}_{l,k}
- \log\int e^{g-t/\sqrt{n} {\mathcal{A}}_{l,k} }.
\end{equation}
%
%Since the last display is invariant by adding a constant to $
%$\zeta_{lk}$ replaced by $\tilde{\cA_{l,k}}$ too.
In the case $1/2<\alpha\le1$, the cost of replacing $f$ by $\log f$
is controlled
by $\llVert g-g_0-(f-f_0)\rrVert_\infty
= \llVert f_0\rho(f/f_0-1)\rrVert_\infty\lesssim\llVert f/f_0 -
1\rrVert_\infty^2$,
which is bounded on $D_n$ by a constant times $\zeta_n^2$ via Lemma
\ref{lem-inter}. Using Remark~\ref{rem-sma} below, the bias (\ref
{sp-bias}) leads to an extra term $\exp\{t \sqrt{n}n^{(-3\alpha
/2)/(1+2\alpha)}\}$ in (\ref{te-lap1}).
%$\left\Vert\right\Vert_{\infty}

%s4.3.2 #&#
\subsubsection{``Uniform'' approximations of efficient influence
functions \texorpdfstring{$\tilde\zeta_{l,k}$}{zeta l,k}}

For any $l\le L_n$ and $k$ between $0$ and $2^{l}-1$, define $
{\mathcal{A}}_{l,k}$ to be the $L^2$-projection of $\zeta_{l,k}$ on
the space spanned by the first $L_n$ levels of wavelet coefficients,
%
%e27 #&#
\begin{eqnarray}
\label{approx1} %\zeta_{n,l,k}^{S}
{\mathcal{A}}_{l,k} &= & \sum
_{1\le\lambda\le L_n} \sum_{0 \le
\mu\le2^{\lambda}-1 } {\langle}
\zeta_{l,k}, \psi_{\lambda\mu} {\rangle}_2 \psi
_{\lambda\mu}. %& - P_{f_0}\left[ \sum_{1\le\la\le L_n} \sum_{0 \le\mu
\end{eqnarray}
For any $l,k$ in the previous ranges, we also set
\[
\tilde{\mathcal{A}}_{l,k} = {\mathcal{A}}_{l,k} -
P_{f_0} {\mathcal{A}}_{l,k}.
\]
%
%le2 #&#
\begin{lem} \label{lem-approx}
Let $f_0$ belong to $ {\mathcal{F}}_0\cap{\mathcal
{C}}^{\alpha}[0,1]$, with $\alpha\ge1$.
For any $l$ such that $1 \le2^l \le2^{L_n}$ and $0\le k\le2^l - 1$,
any density $f$ in $ {\mathcal{F}}$, and $ {\mathcal{A}}_{l,k}$
as in (\ref{approx1}),
\begin{eqnarray*}
\llVert{\mathcal{A}}_{l,k} - \zeta_{l,k} \rrVert
_{\infty} &\lesssim&2^{l(1/2 + \alpha)} 2^{-\alpha L_n},
\\
\biggl\llvert\int_0^1 ( {
\mathcal{A}}_{l,k} - \zeta_{l,k}) (f-f_0) \biggr
\rrvert&\lesssim&\bigl(2^{(l-L_n)\alpha} \wedge2^{-l}\bigr) \llVert
f-f_0\rrVert_2.
\end{eqnarray*}
%
%Let $f_0$ belong to $\cC^{\al}[0,1]$, with $\al\ge1$.
%For any $l$ such that $2^{L_n/2} \le2^l \le2^{L_n}$ and $0\le k\le
%2^l - 1$,
% and any density $f$ in $\cF$,
%For any $l$ such that $2^l\le2^{L_n/2}$ and $0\le k\le2^l - 1$, and
%any density $f$ in $\cF$,
%2^{-\al L_n}, \qquad\left|\int_0^1 (\cA_{l,k} - \zeta_{l,k})(f-f_0)
%2^{-\al L_n/2}\left\Vert f-f_0\right\Vert_2. \]
\end{lem}
\begin{pf}
For any admissible indexes $\lambda,\mu$, let $S_{\lambda\mu}$
denote the support of the wavelet\vspace*{1pt} $\psi_{\lambda\mu}$ in
$[0,1]$ and $\llvert S_{\lambda\mu}\rrvert$ its Lebesgue measure.
The following identity holds both in $L^2[0,1]$ (definition of the
$L^2$-projection) and in $L^\infty[0,1]$ (because $\zeta_{l,k}\in
B_{\infty,\infty}^S$ with $S>0$)
% -is this the minimal condition (?)-)
%
%e28 #&#
\begin{equation}
\label{diff} \zeta_{l,k} - {\mathcal{A}}_{l,k} = \sum
_{\lambda>L_n} \sum_{\mu=0}^{2^\lambda-1}
{\langle}\zeta_{l,k}, \psi_{\lambda\mu} {\rangle}_2
\psi_{\lambda\mu}.
\end{equation}
Since $\psi_{lk}$ belongs to $B_{\infty,\infty}^\alpha$ and $f_0$
to $ {\mathcal{F}}_0\cap{\mathcal{C}}^\alpha$, Lemma~\ref
{lem-hoelder} implies that
$\zeta_{l,k}=\psi_{lk}\cdot f_0^{-1}$ belongs to $B_{\infty,\infty
}^{\alpha}$,
with $\llVert\cdot\rrVert_{\infty,\infty,\alpha}$-norm bounded
above by a
constant times
$\llVert\psi_{lk}\rrVert_{\infty,\infty,\alpha} \llVert
f_0^{-1}\rrVert_{\infty,\infty,\alpha}\lesssim2^{l(1/2+\alpha
)}$, again by
Lemma~\ref{lem-hoelder}, using $f_0^{-1}\in{\mathcal{C}}^\alpha
\subset B_{\infty,\infty}^\alpha$. Now using the localisation
property of the wavelet basis,
\begin{eqnarray*}
\llVert{\mathcal{A}}_{l,k} - \zeta_{l,k} \rrVert
_{\infty} & \le&\Biggl\llVert\sum_{\lambda>L_n} \sum
_{\mu=0}^{2^\lambda-1} {\langle}\zeta_{l,k},
\psi_{\lambda\mu} {\rangle}_2 \psi_{\lambda\mu} \Biggr\rrVert
_{\infty}
\\
& \le&\sum_{\lambda>L_n} 2^{\lambda/2} 2^{-\lambda(\alpha+1/2)}
\max_{0\le\mu\le2^\lambda-1} \bigl[ 2^{\lambda(\alpha+1/2)} \bigl
\llvert{\langle}
\zeta_{l,k}, \psi_{\lambda\mu} {\rangle}_2\bigr\rrvert
\bigr]
\\
& \le&\llVert\zeta_{l,k} \rrVert_{\infty,\infty,\alpha} \sum
_{\lambda>L_n} 2^{-\lambda\alpha} \lesssim2^{l(1/2 + \alpha)}
2^{-\alpha L_n}.
\end{eqnarray*}

Now let us prove that the support of $\zeta_{l,k}- {\mathcal
{A}}_{l,k}$ has diameter at most a constant
times $\llvert S_{lk}\rrvert$. Indeed, $\zeta_{l,k}- {\mathcal{A}}_{l,k}$
written above is a linear combination
of ``high''-frequency wavelets ($\lambda> L_n$), with support diameter
thus at most of the order of
$\llvert S_{\lambda\mu}\rrvert\leq R \llvert S_{lk}\rrvert$, for a
fixed constant $R$, for any
$\lambda>L_n$, any admissible $\mu$, since \mbox{$\lambda>l$}. But in the
sum (\ref{diff}), one may keep only those $\psi_{\lambda\mu}$ whose
support intersects the one of $\zeta_{l,k}$, otherwise the coefficient
${\langle}\zeta_{l,k}, \psi_{\lambda\mu} {\rangle}_2$ is $0$. So,
all supports of the $\psi_{\lambda\mu}$'s which have a nonzero
contribution to (\ref{diff}) are contained in an interval of $[0,1]$
of diameter at most $(2R+1)\llvert S_{lk}\rrvert$. Thus, the diameter
of the support
of $\zeta_{l,k}- {\mathcal{A}}_{l,k}$ is at most $(2R+1)\llvert
S_{lk}\rrvert$.

Now we focus on $\int_0^1 ( {\mathcal{A}}_{l,k}-\zeta
_{l,k})(f-f_0)=\int_0^1 ( {\mathcal{A}}_{l,k}-\zeta
_{l,k})(f-f_0)\mathbh{1}_{\Delta_{l,k}}$, where $\Delta_{l,k}$
denotes the
support of $ {\mathcal{A}}_{l,k}-\zeta_{l,k}$.
%The last integral can also be written$$.
%Cauchy-Schwarz inequality implies
% \[ \int_0^1 (\cA_{l,k}-\zeta_{l,k})(f-f_0) \1_{\Delta_{l,k}}\leq\left
% This implies that there exists a constant $C$ such that the diameter
% of the support of $\zeta_{l,k}-\cA_{l,k}$ is at most $C\left|S_{lk}
%Conclude that
Bounding $ {\mathcal{A}}_{l,k}-\zeta_{l,k}$ by its supremum and
next applying Cauchy--Schwarz inequality,
\begin{eqnarray*}
\biggl\llvert\int_0^1 ( {
\mathcal{A}}_{l,k} - \zeta_{l,k}) (f-f_0)\biggr
\rrvert&\lesssim&\llVert{\mathcal{A}}_{l,k} - \zeta_{l,k}
\rrVert_{\infty} \sqrt{\llvert\Delta_{lk}\rrvert} \llVert
f-f_0\rrVert_2
\\
&\lesssim&2^{(l-L_n)\alpha}\llVert
f-f_0\rrVert_2.
\end{eqnarray*}
To obtain the other part of the bound, the idea is to use a different
approximating sequence $ {\mathcal{D}}_{l,k}$ for which
the comparison to $\zeta_{l,k}$ is easier for large $l$'s.
Define $ {\mathcal{D}}_{l,k}$ to be the function obtained by
replacing $f_0$ in (\ref{def-zet}) by its average
on the support of~$\psi_{lk}$,
%
%e29 #&#
\begin{equation}
\label{approx2} {\mathcal{D}}_{l,k} = \frac{\psi_{lk}}{ [\bar f_0]_{lk}
}, %- P_{f_0}\left[ \frac{\psil}{ [\bar f_0]_{lk}} \right],
\end{equation}
where we have set
\[
[\bar f_0]_{lk} = \frac{1}{\llvert S_{lk}\rrvert} \int
_{S_{lk}} f_0.
\]
Note that since $l\le L_n$, by definition $ {\mathcal{D}}_{l,k}$
belongs to the vector space generated by the first $L_n$ levels of
wavelet coefficients. In particular, it holds $\llVert{\mathcal
{A}}_{l,k}-\zeta_{l,k}\rrVert_2 \le\llVert{\mathcal{D}}_{l,k}-\zeta
_{l,k}\rrVert_2$ by definition of the $L^2$-projection. Since by definition
again $ {\mathcal{D}}_{l,k}-\zeta_{l,k} $ has support included in
$S_{lk}$, one gets
\[
\llVert{\mathcal{A}}_{l,k}-\zeta_{l,k} \rrVert
_2^2 \le\int\mathbh{1}_{S_{lk}} ( {
\mathcal{D}}_{l,k}-\zeta_{l,k})^2 \le\llVert{
\mathcal{D}}_{l,k}-\zeta_{l,k}\rrVert_{\infty}^2
\llvert S_{lk}\rrvert.
\]
Next, one bounds the last sup-norm. Denoting $\rho_0:= \inf_{x\in
[0,1]} f_0(x)$,
\begin{eqnarray*}
\llVert{\mathcal{D}}_{l,k} - \zeta_{l,k} \rrVert
_{\infty} & \le&\rho_0^{-1} \llVert
\psi_{lk}\rrVert_{\infty} \sup_{x\in S_{lk}} \bigl
\llvert f_0(x)- [\bar f_0]_{lk}\bigr\rrvert
\\
& \le&\rho_0^{-1} \llVert\psi_{lk}\rrVert
_{\infty} \sup_{x\in S_{lk}} \llvert S_{lk}\rrvert
^{-1} \biggl\llvert\int_{S_{lk}}
\bigl(f_0(x)- f_0(u)\bigr)\,du\biggr\rrvert
\\
& \le&\rho_0^{-1} \llVert\psi_{lk}\rrVert
_{\infty} \sup_{x\in S_{lk}} \llvert S_{lk}\rrvert
^{-1} \int_{S_{lk}} \llvert x-u\rrvert \,du
\\
& \le&\rho_0^{-1} \llVert\psi_{lk}\rrVert
_{\infty} \llvert S_{lk}\rrvert^{-1} \llvert
S_{lk}\rrvert^{2}/2 \lesssim2^{-l/2},
\end{eqnarray*}
where for the third inequality we have used the fact that $f_0$ is at
least H\"{o}lder $1$ and for the last inequality that $\llvert
S_{lk}\rrvert$ is of
the order $2^{-l}$ and $\llVert\psi_{lk}\rrVert_{\infty}$ of the order
$2^{l/2}$ up to constants.
Thus,
%Since, as noted above, $\cA_{l,k} - \zeta_{l,k}$ has support $
%$\left|\Delta_{lk}\right| \leqa\left|S_{lk}\right|\leqa2^{-l}$,
%
\begin{eqnarray*}
\biggl\llvert\int_0^1 ( {
\mathcal{A}}_{l,k} - \zeta_{l,k}) (f-f_0)\biggr
\rrvert& \le&\llVert{\mathcal{A}}_{l,k} - \zeta_{l,k} \rrVert
_{2} \llVert f-f_0 \rrVert_2
\\
& \le&\sqrt{\llvert S_{lk}\rrvert} \llVert{\mathcal{D}}_{l,k}
- \zeta_{l,k} \rrVert_{\infty} \llVert f-f_0
\rrVert_2
\\
&\lesssim& 2^{-l} \llVert f-f_0\rrVert
_2.
\end{eqnarray*}\upqed
\end{pf}
%
%re1 #&#
\begin{remark} \label{rem-sma}
In\vspace*{1pt} the case $1/2<\alpha< 1$, similarly one gets the bound
$(2^{-l\alpha} \wedge(2^{l-L_n})^\alpha)\llVert f-f_0\rrVert_2$.
The minimum of the bounds is attained for $2^l=2^{L_n/2}$.
%Then using the two previous approximations in the regimes $2^l>M_n$
%and $2^l\le M_n$ respectively, one
This leads to a bound for the integral $\int_0^1 ( {\mathcal
{A}}_{l,k}-\zeta_{l,k})(f-f_0)$
equal to\break  $n^{(-3\alpha/2)/(1+2\alpha)}$ for all considered indexes $l,k$.
\end{remark}

%s4.3.3 #&#
\subsubsection{Change of variables}

Now everything is in place to start exploiting \mbox{(\ref{te-lap1})--(\ref{te-gt1})}. First, rewrite (\ref{te-lap1}) as
%
%e30 #&#
\begin{equation}
\label{cha} \frac{\int e^{ \ell_n(f_t) - \ell_n(f_0) } \,d\Pi^{D_n}(f) }{
\int e^{ \ell_n(f) - \ell_n(f_0) } \,d\Pi^{D_n}(f) } = \frac{\int\mathbh
{1}_{D_n}(f) e^{ \ell_n(f_t) - \ell_n(f_0) } \,d\Pi
(f) }{
\Pi(D_n\given X^{(n)}) \int e^{ \ell_n(f) - \ell_n(f_0) } \,d\Pi(f) }.
\end{equation}
The expression $\log f_t=(\ref{te-gt1})$, due to its invariance by
adding a constant and recalling that $g=T-c(T)$ from (\ref
{def-prior-log}), can be seen as a function of $T-t {\mathcal
{A}}_{l,k}/\sqrt{n}$ [the constant $c(T)$ vanishes]. More precisely,
we are now ready to change variables in the prior by setting
%
%e31 #&#
\begin{equation}
\label{def-G} \tilde T = T-\frac{t}{\sqrt{n}} {\mathcal{A}}_{l,k}.
\end{equation}
Essentially, if the ``complexity''
of $ {\mathcal{A}}_{l,k}$ is not too large in view of the chosen
prior~$\Pi$, the fact of having $f_t$ instead of $f$ in (\ref
{te-lap1}) will not matter much and the corresponding ratio of
integrals will be close to $1$. We treat the case of log-Lipschitz
priors on coefficients first. In fact, as can intuitively be guessed, a
prior with heavy tails is less influenced under shift transformations
than a more concentrated prior.

Denote by $ {\mathcal{C}}_n=\{(\lambda,\mu)\in\mathbb{N}^2, 0\le\mu\le
2^{\lambda}-1, 1\le2^\lambda\le2^{L_n}\}$. By the definition~(\ref
{approx1}) of $ {\mathcal{A}}_{l,k}$,
\begin{eqnarray*}
{\langle} {\mathcal{A}}_{l,k},\psi_{\lambda\mu} {
\rangle}_2 & =& {\langle}\zeta_{l,k},\psi_{\lambda\mu} {
\rangle}_2, \qquad(\lambda,\mu)\in{\mathcal{C}}_n, (l,k)\in{\mathcal
{C}}_n. %\psg\cA_{l,k},\psi_{\la\mu} \psd_2 & =& \delta_{l\la}\delta_{k
\end{eqnarray*}

\paragraph*{Log-Lipschitz prior} With the chosen prior on $f$, the numerator
in (\ref{cha}) is in fact an integral over the law of the coefficients
of $T$ in (\ref{def-T}),
that is, over (a subset of)~$\mathbb{R}^{2^{L_n}}$.
The change of variables (\ref{def-G}) is thus a shift in $\mathbb
{R}^{2^{L_n} }$, and its Jacobian is~$1$. The coordinates of $T$ in the
wavelet basis $\{\psi_{\lambda\mu}\}$ have densities $\sigma
_\lambda^{-1}\varphi(\theta_{\lambda,\mu}/\sigma_\lambda)$ with
respect to $d\theta_{\lambda,\mu}$ (we denote by $\theta_{\lambda,\mu}$
the integrating variable). The transformation in density can be
controlled by, since $\varphi=\varphi_H$ has a Lipschitz logarithm,
\begin{eqnarray*}
&& \log\prod_{(\lambda,\mu)\in{\mathcal{C}}_n} \frac{\varphi
(\theta_{\lambda,\mu}/\sigma_\lambda)}{
\varphi(\{\theta_{\lambda,\mu}-t/{\sqrt{n}}\langle
{\mathcal{A}}_{l,k},\psi_{\lambda\mu}{\rangle}_2 \}/\sigma
_\lambda)}
\\
&&\qquad  =
\sum_{(\lambda,\mu)\in{\mathcal{C}}_n} (\log\varphi) (\theta
_{\lambda,\mu}/\sigma_\lambda) - (\log\varphi) \biggl(\biggl\{
\theta_{\lambda,\mu}-\frac{t}{\sqrt
{n}}{\langle} {\mathcal{A}}_{l,k},
\psi_{\lambda\mu}{\rangle}_2\biggr\}\Big/\sigma_\lambda\biggr)
\\
&&\qquad  \le\sum_{(\lambda,\mu)\in{\mathcal{C}}_n} \frac{\llvert t\rrvert
}{\sqrt
{n}\sigma_\lambda}
\bigl\llvert{\langle}\zeta_{l,k},\psi_{\lambda\mu} {
\rangle}_2 \bigr\rrvert\qquad \forall(l,k)\in{
\mathcal{C}}_n.
\end{eqnarray*}
We now study conditions on the $\sigma_\lambda$'s under which the
last display is bounded above by $C\llvert t\rrvert$.
Let us split the sum over $ {\mathcal{C}}_n$ in the two cases
$\lambda\le l$ and $\lambda>l$. When $\lambda\le l$, for any fixed
level $\lambda$ there is a constant number of wavelets $\psi_{\lambda
\mu}$ intersecting the support of $\zeta_{l,k}$. Combined with
$\llvert{\langle}\zeta_{l,k},\psi_{\lambda\mu}{\rangle}_2\rrvert\leq
\rho
_0^{-1}$, this leads to the condition $\sum_{\lambda\le l} \sigma
_{\lambda}^{-1} \lesssim\sqrt{n}$. When $\lambda>l$, for any fixed
level $\lambda$ there is a constant times $2^{\lambda-l}$ wavelets
$\psi_{\lambda\mu}$ intersecting the support of $\zeta_{l,k}$,
leading to
\begin{eqnarray*}
\sum_{(\lambda,\mu) \in{\mathcal{C}}_n} \sigma_\lambda^{-1}
\bigl\llvert{\langle}\zeta_{l,k},\psi_{\lambda\mu} {
\rangle}_2 \bigr\rrvert& \le&\sum_{l<\lambda\le L_n}
\sigma_\lambda^{-1} 2^{\lambda-l} 2^{-\lambda(1/2+\alpha)} \llVert
\zeta_{l,k}\rrVert_{\infty,\infty,\alpha}
\\
& \le&\sum_{l<\lambda\le L_n} \sigma_\lambda^{-1}
2^{(l-\lambda
)(\alpha-1/2)},
\end{eqnarray*}
where we have used that $\zeta_{l,k}=\psi_{lk}/f_0$ is $\alpha
$-smooth and applied
Lemma~\ref{lem-hoelder}. These conditions are quite mild. In
particular, for $\alpha>1/2$ they are implied by $\sum_{\lambda\le
L_n}\sigma_\lambda^{-1} \lesssim\sqrt{n}$.

\paragraph*{Gaussian prior} %Here we have to restrict slightly the range of
%possible $\sigma_l$.
Let us write explicitly the log-ratio of densities
%[See if previous bounds give the same condition as the implicit
%Gaussian calculus, using the fact that $f$ and $f_0$ are close. Since
%we work directly on series priors, one can make explicit bounds at the
%coefficients level.]
%
%e32 #&#
\begin{eqnarray}\label{chg}
&&\sum_{(\lambda,\mu)\in{\mathcal{C}}_n} (\log\varphi) (\theta
_{\lambda,\mu}/\sigma_\lambda) - (\log\varphi) \biggl(\biggl\{
\theta_{\lambda,\mu}-\frac{t}{\sqrt
{n}}{\langle} {\mathcal{A}}_{l,k},
\psi_{\lambda\mu}{\rangle}_2\biggr\}\Big/\sigma_\lambda\biggr)
\nonumber\\[-8pt]\\[-8pt]\nonumber
&&\qquad  =\sum_{(\lambda,\mu)\in{\mathcal{C}}_n} \frac{t^2}{2n\sigma
_\lambda^2} {
\langle} {\mathcal{A}}_{l,k},\psi_{\lambda\mu
}{\rangle}_2^2
-\frac{t}{\sqrt{n}\sigma_\lambda^2}\theta_{\lambda,\mu} {\langle}
{\mathcal{A}}_{l,k},
\psi_{\lambda\mu}{\rangle}_2.
\end{eqnarray}
The obtained quantity still depends on the integrating variables
$\theta_{\lambda,\mu}$. The idea is to exploit the fact that on
$D_n$, it holds $\llVert g-g_0\rrVert_2\lesssim\varepsilon_n$, which
is obtained
along the way in the proof of Lemma~\ref{lem-inter}. But
$g=T-c(T)=T-c(T)1$, where $1$ denotes the constant function equal to
$1$. Since $1$ is orthogonal to high levels of wavelet coefficients, it
means that for large enough $\lambda$, say $\lambda> K$, and any $\mu
$, it holds
$\llvert\theta_{\lambda,\mu}-g_{0,\lambda,\mu}\rrvert\le\llVert
g-g_0\rrVert_2\lesssim
\varepsilon_n$. So, for such $\lambda,\mu$, we decompose
$\theta_{\lambda,\mu} = \theta_{\lambda,\mu}-g_{0,\lambda,\mu
}+g_{0,\lambda,\mu}$. For coefficients
$\theta_{\lambda,\mu}$ such that $\lambda\le K$, we use a different
argument.

Let us first deal with the term containing $\theta_{\lambda,\mu}$ in
(\ref{chg}) when $\lambda\le K$. From the beginning of the proof,
using Lemma~\ref{lem-max}, one can restrict slightly the set $D_n$ by
intersecting it with the set
$\{T\dvtx  \max_{\lambda\le K,\mu} \llvert{\langle}T,\psi_{lk}{\rangle
}\rrvert
\le C\sqrt{n}\varepsilon_n\}$.
%This implies that $\max_{\la\le K,\mu}\left|\te_{\la,k}\right|\leqa\rn
Hence, using that $\llvert{\langle} {\mathcal{A}}_{l,k},\psi_{\lambda
\mu}{\rangle}_2\rrvert\lesssim1$ and the assumed specific form of
$\sigma
_\lambda$, one gets that the term at stake is at most a fixed constant
times $\llvert t\rrvert\varepsilon_n$.
%To handle this dependence, we write $g_{\la,\mu} = g_{\la,\mu}-g_{0,
%the first difference is bounded by $\left\Vert g-g_0\right\Vert_2=
%which itself is bounded by a constant times $\left\Vert f-f_0\right
%using the fact that we work on a set where the sup-norm distance
%already tends to $0$.

Now we bound (\ref{chg}) and only have to deal with $\lambda>K$ for
the part depending on $\theta_{\lambda,\mu}$.
%When $(l,k)\in\cC_n^\iota$, the first term in \eqref{chg} leads again
%to the condition $\rn\sigma_l\ge1$. For the second term, using the
%previous decomposition inserting $g_0$ and $\left|\te_{\la,\mu}-g_{0,
For any $(l,k)\in{\mathcal{C}}_n$, the first term in (\ref{chg}) is
\[
\frac{t^2}{2n}\sum_{(\lambda,\mu)\in{\mathcal{C}}_n} \sigma
_\lambda^{-2} {\biggl\langle}\frac{\psi_{lk}}{f_0},
\psi_{\lambda\mu
}{\biggr\rangle}_2^2 \lesssim
\frac{t^2}{n}\sum_{\lambda\le l}\sigma_\lambda^{-2}
+ \frac{t^2}{n}\sum_{l<\lambda\le L_n} \sigma_\lambda^{-2}
2^{(l-\lambda)2\alpha},
\]
where the bound is obtained in a similar way as in the log-Lipschitz
case by distinguishing the cases
$\lambda\le l$ and $l>\lambda$. For the second term, we decompose
$\theta=\theta-g_0+g_0$ and use Cauchy--Schwarz inequality on the
$\theta-g_0$ part,% to make appear
\begin{eqnarray*}
&& \frac{\llvert t\rrvert}{\sqrt{n}}\llVert g-g_0\rrVert_2\biggl\{ \sum
_{ (\lambda,\mu)\in
{\mathcal{C}}_n} \sigma_\lambda^{-4} {\biggl
\langle}\frac{\psi
_{lk}}{f_0},\psi_{\lambda\mu}{\biggr\rangle}_2^2
\biggr\}^{1/2}
\\
&&\qquad \le\llvert t\rrvert\frac
{\varepsilon_n}{\sqrt{n}} \biggl\{ \sum
_{ \lambda\le l} \sigma_\lambda^{-4} + \sum
_{l<\lambda\le
L_n} \sigma_\lambda^{-4}
2^{(l-\lambda)2\alpha} \biggr\}^{1/2}.
\end{eqnarray*}
Finally, the remaining term with $g_0$ is bounded by $\llvert t\rrvert
/\sqrt{n}$ times
\begin{eqnarray*}
&& \sum_{ (\lambda,\mu)\in{\mathcal{C}}_n} \sigma_\lambda^{-2}
\biggl\llvert g_{0,\lambda\mu} {\biggl\langle}\frac{\psi_{lk}}{f_0},
\psi_{\lambda\mu}{\biggr\rangle}_2\biggr\rrvert
\\
&&\qquad \le\sum
_{ \lambda\le l} \sigma_\lambda^{-2}2^{-\lambda(1/2+\alpha)}
+ \sum_{l<\lambda\le L_n} \sigma_\lambda^{-2}
2^{(l-\lambda
)2\alpha}2^{-\lambda(1/2+\alpha)}.
\end{eqnarray*}
Under the condition of the theorem $\sigma_\lambda\ge2^{-\lambda
(1/4+\alpha)}$ for any $0\le\lambda\le L_n$, all bounds
obtained above in the Gaussian case are less than $C(\llvert t\rrvert+t^2)$.

%s4.3.4 #&#
\subsubsection{End of proof}

To conclude for both considered classes of priors, %note that the
%presence of the integrating set $D_n$ does not matter anymore once the
%ratios of densities above are controlled.
note that the indicator $\mathbh{1}_{D_n}(f)$ in (\ref{cha}) becomes
$\mathbh{1}
_{D_n'}(f_t)$ under the change of variables---for a set $D_n'$ that
one can write explicitly, although this will not be needed here---and
one simply further bounds the indicator $\mathbh{1}_{D_n'}$ by $1$ on the
numerator. Once the change of variable is done, the assumed conditions
on $\{\sigma_{l}\}$ ensure that the ratio of densities is bounded by
$e^{C(\llvert t\rrvert+t^2)}$ for some constant $C$ and one gets
\[
\frac{\int e^{ \ell_n(f_t) - \ell_n(f_0) } \mathbh{1}_{D_n}(f) \,d\Pi
(f) }{
\int e^{ \ell_n(f) - \ell_n(f_0) } \,d\Pi(f) } \le e^{C(\llvert t\rrvert
+t^2)} \frac{\int e^{ \ell_n(f) - \ell_n(f_0) } \,d\Pi
(f) }{
\int e^{ \ell_n(f) - \ell_n(f_0) } \,d\Pi(f) } = e^{C(\llvert t\rrvert+t^2)}.
\]
Inequality (\ref{te-lap1}) can thus be further written, again for a
fixed constant $C$, %pushing the term in $W_n(\cA_{l,k})$ to the other
%side of the inequality,
%
\[
%e^{-t W_n(\cA_{l,k}) } E^{\Pi^{D_n}}\left[ e^{t\rn\psg g-g_0,
E^{\Pi^{D_n}} \bigl[ e^{t\sqrt{n}{\langle}g-\Gamma^{L_n},\psi_{lk}
{\rangle}_2} \given
X^{(n)} \bigr] \le e^{C(\llvert t\rrvert+t^2)} \Pi\bigl(D_n\given
X^{(n)}\bigr)^{-1}.
\]
%
%Control the change in the indexing set. The integrating set $D_n$ is
%also altered via the change of variables. To control this change, note
%that $D_n$ verifies, using the triangle inequality,
% \left\Vert f_t-f_0\right\Vert_2 \le\zeta_n+ \left\Vert f_t-f\right
%Thanks to Lemma~\ref{lem-cont}, the terms $\left\Vert f_t-f\right\Vert
%_{\infty}$ and $
%by
%Now it remains to control the ratio of integrals in \eqref{te-lap1}
%once the change of variables has been done. The quantity $e^{
%$e^{\ell_n(f) - \ell_n(f_0)}$. We are left with
For any $s>0$, similar to the Gaussian white noise case,
\begin{eqnarray*}
\hspace*{-3pt}&& \int\bigl\llVert g^{L_n} - \Gamma^{L_n} \bigr\rrVert
_{\infty} \,d\Pi^{D_n}\bigl(f\given X^{(n)}\bigr)
\\
\hspace*{-3pt}&&\quad \le\sum
_{l\le L_n} \frac{2^{l/2}}{\sqrt{n}} E^{\Pi^{D_n}}\Bigl[
\max_{0 \le k \le2^{l}-1} \bigl\llvert\sqrt{n} {\bigl\langle}g-
\Gamma^{L_{n}},\psi_{lk} {\bigr\rangle}_2\bigr
\rrvert\big\given X^{(n)}\Bigr]
\\
\hspace*{-3pt}&&\quad  \le \sum_{l\le L_n} 2^{l/2}
\frac{1}{\sqrt{n}s} \log\Biggl[ \sum_{k=0}^{2^{l}-1}
E^{\Pi^{D_n}}\bigl[ e^{s\sqrt{n}{\langle}g-\Gamma^{L_n},\psi_{lk}
{\rangle}_2} + e^{-s\sqrt{n}{\langle}g-\Gamma^{L_n},\psi_{lk} {\rangle
}_2} \given
X^{(n)}\bigr] \Biggr]
\\
\hspace*{-3pt}&&\quad \lesssim \sum_{l\le L_n} 2^{l/2}
\frac{1}{\sqrt{n}s} \log\bigl[ 2^l e^{C(s+s^2)}\bigr] + \sum
_{l\le L_n} 2^{l/2} \frac{1}{\sqrt{n}s} \log
\frac{1}{\Pi
(D_n\given X^{(n)})}. % \leqa\frac{1}{\rn}
% \sum_{l\dvtx  2^l\le J_n} \sqrt{l} 2^{l/2} \leqa\veps_{n,\alpha}^*,
\end{eqnarray*}
Set $s=\sqrt{l}$. The first term in the last display is bounded by a
constant times $\frac{1}{\sqrt{n}}
\sum_{l\le L_n} \sqrt{l} 2^{l/2} \lesssim\varepsilon_{n,\alpha
}^*$. Now coming back to the application of Markov's inequality one
gets, with $\Xi$ the function $\Xi\dvtx u\to u\log u^{-1}$,
\[
E_{f_0}^n \Pi\bigl[\llVert g-g_0\rrVert
_\infty>M_n \varepsilon_{n,\alpha}^* \given
X^{(n)}\bigr] \lesssim M_n^{-1} +
M_n^{-1}E_{f_0}^n \Xi\bigl(\Pi
\bigl(D_n \given X^{(n)}\bigr)\bigr)+ o(1).
\]
With $M_n\to\infty$, the fact that $\Xi$ is bounded on $[0,1]$
completes the proof.
%?????????????????? \qed

In the case $1/2<\alpha\le1$, the only difference is that one gets
two extra terms: one from going at the logarithmic level, which
eventually leads to a rate $\zeta_n^2$; another one from the
semiparametric bias
in (\ref{sp-bias}), which leads to a rate $\rho_n=n^{(1/2-3\alpha/2)/(1+2\alpha)}$. This leads to a sup-norm rate of $\zeta
_n^2\vee\rho_n=\zeta_n^2$.

%In the case $1/2<\al\le1$,
Once the rate $\zeta_n^2$ has been obtained, one can restart the proof
once again, but this time knowing that one can use a better
intermediate rate of $\zeta_n^2$ in sup-norm.
One can then write
\[
\biggl\llvert\int\rho\biggl(\frac{f-f_0}{f_0}\biggr)\psi_{lk}
\biggr\rrvert\le\biggl\llVert\rho\biggl(\frac
{f-f_0}{f_0}\biggr)\biggr\rrVert
_\infty\llVert\psi_{lk}\rrVert_1 \lesssim\biggl
\llVert\frac{f-f_0}{f_0} \biggr\rrVert_\infty^2
2^{-l/2} \lesssim\bigl(\zeta_n^2
\bigr)^2 2^{-l/2}.
\]
This eventually leads to the accelerated rate $(\zeta_n^2)^2\vee\rho
_n$. Iterating
this procedure leads to $(\zeta_n^2)^{2^p}\vee\rho_n$, any $p\ge1$
and, for any given $\alpha>1/2$, to $\rho_n$ as
final rate, up to logarithmic terms.

%s4.4 #&#
\subsection{Density estimation, dyadic histogram priors}

We follow the scheme of proof used for uniform priors in white noise,
this time in density estimation, with the wavelet basis of expansion
being the Haar system. The very specific properties of the Haar basis,
particularly its close links to approximation by dyadic histograms,
enable a simplified argument. In particular, as we demonstrate below,
the semiparametric bias is always negligible, provided the parameters
of the prior are reasonably chosen.

Let us set
\[
\hat{f}_{lk} = \mathbb{P}_n \psi_{lk}^H
= \frac{1}{n} \sum_{i=1}^n
\psi_{lk}^H(X_i).
\]
Set $D_n=\{f, h(f,f_0)\le\varepsilon_n\}$, where here $\varepsilon
_n=\varepsilon_{n,\alpha}^*$ up to multiplication by a~large enough
constant. Lemma~\ref{lem-histr} implies
that $E_{f_0}^n\Pi[D_n\given X^{(n)}]\to1$.
Let $h_n, L_n$ be defined as in (\ref{def-jn}), and denote by
$f^{L_n}$ the projection of $f$ onto
the subspace $ {\mathcal{V}}_{L_n}:=\operatorname{Vect}\{\varphi
^H, \psi_{lk}^H, l < L_n, 0\le k<2^l\}$ and $f^{L_n^c}$ the
projection of $f$ onto
$\operatorname{Vect}\{\psi_{lk}^H, l \ge L_n, 0\le k<2^l\}$ (i.e.,
for simplicity we keep the notation $f^{L_n}$ from Section~\ref
{firstex}, although\vspace*{2pt} the basis of projection is now the Haar basis \emph{and} $l<L_n$ replaces $l\le L_n$). If $\hat f^{L_n}$ denotes the
element of $ {\mathcal{V}}_{L_n}$ of coordinates
$\{\hat f_{lk}\}$ in the basis $\{\psi_{lk}^H\}$, one can bound
$E^{\Pi^{D_n}}[ \llVert f-f_0\rrVert_{\infty} \given X^{(n)} ]$ from
above by
%f^{J_n^c} -f_0^{J_n^c}. \]
%As before,
%E^{\Pi^{D_n}}[ \left\Vert f-f_0\right\Vert_{\infty} \given X^{(n)} ]
%
\[
\underbrace{\int\bigl\llVert f^{L_n} -\hat{f}{}^{L_n} \bigr
\rrVert_{\infty} \,d\Pi^{D_n}\bigl(f\given X^{(n)}
\bigr)}_{(\mathrm{i})} + \underbrace{\int\bigl\llVert f^{L_n^c}\bigr
\rrVert_{\infty} \,d\Pi\bigl(f\given X^{(n)}\bigr)}_{(\mathrm{ii})}
+ \underbrace{\bigl\llVert\hat{f}{}^{L_n}-f_0\bigr\rrVert
_{\infty}}_{(\mathrm{iii})}.
\]
The term (iii) can be bounded by $\llVert\hat{f}{}^{L_n} -
f_0^{L_n}\rrVert
_{\infty}+\llVert f_0^{L_n^c}\rrVert_{\infty}$. The second term in
this sum is
pure bias and the first term is
\[
\Biggl\llVert\bigl(\mathbb{P}_n\varphi^H -
P_{f_0} \varphi^H\bigr) \varphi^H(\cdot) + \sum
_{l=0}^{L_n-1}\sum
_{k=0}^{2^l-1} \bigl(\mathbb{P}_n
\psi_{lk}^H - P_{f_0} \psi_{lk}^H
\bigr) \psi_{lk}^H(\cdot) \Biggr\rrVert_{\infty}.
\]
This term is bounded in expectation by $\varepsilon_{n,\alpha}^*$,
exactly as in Lemma~\ref{lem-Gamma}.

Next, the high-frequency bias term (ii) is zero. Indeed, for any draw
$f$ from the prior, in the inner-product ${\langle}f, \psi_{lk}^H
{\rangle}_2$, the first element is a dyadic
histogram at resolution level $L_n$, so is constant over the support of
the Haar basis element
$\psi_{lk}^H$ if $l\ge L_n$. Hence, the previous inner-product is zero
$\Pi$-almost surely, and thus also $\Pi[\cdot\given X^{(n)}]$ almost surely.

Second, one studies ${\langle}f^{L_n} - \hat f^{L_n}, \psi_{lk}^H
{\rangle}_2$ in the BvM-regime
$l < L_n$.
Following the maximal inequality approach from Section~\ref{firstex},
it is enough to bound the
posterior expectation of $\exp(t\sqrt{n}{\langle}f^{L_n}-\hat
f^{L_n}, \psi^H_{lk}{\rangle}_2)$, for any possible $k$ and $l<L_n$
and say $\llvert t\rrvert\lesssim\log n$ (also,\vspace*{1pt} to simplify the
notation below we
omit mentioning the scaling function $\varphi^H$, but the same Laplace
transform control is obtained for it in a similar way).
To do so, apply Lemma~\ref{lem-lapu} with $\gamma_n=\psi_{lk}^H$,
for any given $l,k$ with
$l < L_n$.
The conditions of the lemma are satisfied with $a_n=\varepsilon_n$,
since $\psi_{lk}^H$ is bounded in $L^2[0,1]$ and has a sup-norm
bounded by a constant times $2^{l/2}\lesssim2^{L_n/2}=o(1/\varepsilon
_n)$ if $\alpha>1/2$. %(if $\al=1/2$ one adapts slightly the argument
%by taking a smaller `$t$' in the final balancing of the maximal
%inequality, we omit the details).
Noting that, again for $l,k$ with $l < L_n$,
\begin{eqnarray*}
{\bigl\langle}f^{L_n}-\hat f^{L_n}, \psi^H_{lk}
{\bigr\rangle}_2 &=& {\bigl\langle}f^{L_n}-f^{L_n}_0, \psi^H_{lk} {\bigr
\rangle}_2 + {\bigl
\langle}f^{L_n}_0 -\hat f^{L_n},
\psi^H_{lk} {\bigr\rangle}_2
\\
&=& {\bigl
\langle}f-f_0, \psi^H_{lk} {\bigr
\rangle}_2 - W_n\bigl(\psi_{lk}^H
\bigr),
\end{eqnarray*}
an application of Lemma~\ref{lem-lapu} leads to, with $\Pi^{D_n}$ the
restriction of $\Pi$ to $D_n$,
\[
E^{\Pi^{D_n}}\bigl[ e^{t\sqrt{n}{\langle}f^{L_n}-\hat f^{L_n}, \psi
^H_{lk}{\rangle}_2 } \given X^{(n)} \bigr]
\lesssim e^{ Ct^2 } \frac{ \int e^{\ell_n(f_t)-\ell_n(f_0) } \,d\Pi^{D_n}(f)
}{\int e^{\ell_n(f)-\ell_n(f_0)}\,d\Pi^{D_n}(f) },
\]
where $f_t$ is defined by $\log f_t = \log f - t\psi_{lk}^H/\sqrt{n}-
c(f e^{- t\psi_{lk}^H/\sqrt{n}})$ (again, having $\tilde\psi
_{lk}^H$ or $\psi_{lk}^H$ at both places in the last equality does not
matter since the constant simplifies).
Below to simplify the notation, we denote $\gamma_n=\psi_{lk}^H$.

In the last display, once coming back to $\Pi$ via $d\Pi
^{D_n}(f)=\mathbh{1}
_{D_n}\,d\Pi(f)/\Pi(D_n)$, the variable $f$ is a random dyadic
histogram over the subdivision with intervals $I_{\mu}^{L_n}=(\mu
2^{-L_n},(\mu+1)2^{-L_n})$ and $0\le\mu\le2^{L_n}-1$. Denote by
$\gamma_{n,\mu}$ the value of $\gamma_n$ over the interval $I_\mu^{L_n}$.
Next, observe that both $f_t$ and $f$ are, writing the histogram prior
in terms of its coefficients over the subdivision, functions of $\omega
=(\omega_0,\ldots,\omega_{2^{L_n}-1})$ and the integral over $f$ is
nothing but an integral over $\omega\in{\mathcal{S}}_{L_n}$. On
the other hand, from the expression of $f_t$,
%The value taken by $f_t$ on $I_\mu$ is a function of $\omega_\mu e^{-t
using the fact that $\psi_{lk}^H$ is constant over each individual
interval $I_\mu$ (since $l<L_n$),
one sees that $f_t$ is a dyadic histogram over $I_\mu$ with weights
given by the vector $\zeta$
%
%e33 #&#
\begin{equation}
\label{cvarsim} \zeta:=(\zeta_\mu)_{0\le\mu\le2^{L_n}-1}= \biggl(
\frac{\omega_\mu
e^{-t\gamma_{n,\mu}/\sqrt{n}}}{
\sum_{\mu} \omega_\mu e^{-t\gamma_{n,\mu}/\sqrt{n}}} \biggr)_{0\le
\mu\le2^{L_n}-1}.
\end{equation}
Now we are in position to change variables in the above expression by
taking $\zeta$ as the new variable (this technique is developed in
\cite{cr13} for general, \emph{fixed} influence functions). The change
of variables introduces a multiplicative factor $M(\zeta)$ in front of
$d\Pi(\zeta)$, factor which is the product of the variation in
density of the Dirichlet law under the change of variable
and the Jacobian of the change of variable,
\[
\frac{ \int e^{\ell_n(f_t)-\ell_n(f_0) } \,d\Pi^{D_n}(f)
}{\int e^{\ell_n(f)-\ell_n(f_0)}\,d\Pi^{D_n}(f) } = \Pi\bigl[D_n\given X^{(n)}
\bigr]^{-1} \frac{ \int_{\tilde D_n} e^{\ell
_n(f(\zeta))-\ell_n(f_0) } M(\zeta) \,d {\mathcal{D}}_\alpha
(\zeta)
}{\int e^{\ell_n(f(\omega))-\ell_n(f_0)}\,d {\mathcal{D}}_\alpha
(\omega) },
\]
%
%where in a slight abuse of notation we denote $d\Pi(\omega)$ instead
%of $d\cD_\al(\omega)$ (recall that by definition $\Pi$ is the
%push-forward measure of $D_\al$ under $\omega\to f(\omega)$).
where $\tilde D_n$ is the new integrating set after change of variables
and the notation $f(\zeta)$ is used for
$f(\zeta)(\cdot) =
2^L \sum_{\mu=0}^{2^L-1} \zeta_\mu\mathbh{1}_{I_k^L}(\cdot)$ and similarly
for $f(\omega)$.
Computation of the Jacobian say $\Delta(\zeta)$ is done in Lemma~\ref
{lem-jaco}.
Calculating $M(\zeta)=d {\mathcal{D}}_\alpha(\omega)/d
{\mathcal{D}}_\alpha(\zeta)\Delta(\zeta)$ gives that
%the multiplicative term
$M(\zeta)$ satisfies %(see \cite{cr13} for details in the fixed
%influence function case)
%
\begin{eqnarray*}
M(\zeta) e^{-t/{\sqrt{n}} \sum_{\mu} \alpha_\mu\gamma
_{n,\mu} } &=& \biggl[ \int_0^1
e^{t\gamma_n(x)/\sqrt{n}}f(\zeta) (x)\,dx \biggr]^{-\sum_{\mu} \alpha_\mu
}
\\
&=& \biggl[ \int
_0^1 e^{-t\gamma_n(x)/\sqrt{n}}f(\omega) (x)\,dx
\biggr]^{\sum_{\mu} \alpha_\mu}.
\end{eqnarray*}
For the term on the right-hand side, since $\alpha>1/2$ it holds
$\llvert t\rrvert\llVert
\gamma_n\rrVert_\infty/\sqrt{n}=o(1)$ so one can expand the exponential
function by writing
$e^u=1+O(u)$ as $u=o(1)$. Next, write $f=f_0+(f-f_0)$, so that the expression
under brackets writes $1 +O((t/\sqrt{n})[ {\langle}f_0,\gamma
_n{\rangle}_2 +
{\langle}f-f_0,\gamma_n{\rangle}_2])$. The last term is a
$O(\llvert t\rrvert/\sqrt{n})$, since the Haar-coefficients of $f_0$
are certainly
bounded and those of $f-f_0$ are bounded above by a constant times
$\llVert
\gamma_n\rrVert_{\infty}h(f,f_0)$ which is bounded because
$2^{L_n/2}\varepsilon_n=O(1)$. So if
$\sum_\mu\alpha_\mu/\sqrt{n}=O(1)$, the term at stake is $O(1)$.
But this condition follows from (\ref{prk}), because the number of
terms in the previous sum is $2^{L_n}=O(\sqrt{n})$ when $\alpha>1/2$.
%(take $\ga=\varphi^H$), in order get the bound $e^{\left|t\right|\sqrt{

Now we deal with the exponential term on the left-hand side of the last
display. A term in the sum $\sum_{\mu} \alpha_\mu\gamma_{n,\mu}$
is nonzero only if the support of the Haar basis element $\gamma
_n=\psi_{lk}^{H}$
(or $\varphi^H$) intersects $I_{\mu}^{L_n}$. This is the case for
$2^{L_n-l}$ terms for $\psi_{lk}^H$ and
$2^{L_n}$ terms for $\varphi^H$. Using $\llVert\psi_{lk}^H\rrVert
_\infty
\lesssim2^{l/2}$, this shows\vspace*{1pt} that the considered sum is at most a
constant times $2^{L_n-l/2}\lesssim2^{L_n}$. In particular, for
$\alpha>1/2$ the considered term
is a $O(1)$.

The previous reasoning shows that the change of variable part generates
a multiplicative factor $O(1)$, times the term $e^{Ct^2}$ coming from
the application of Lemma~\ref{lem-lapu}. To conclude, it now suffices
to apply the maximal inequality technique in the same way as for
log-densities prior.
%?????????????????? \qed
%Semiparametric bias term is negligible
%$\rn\int(f_0-f_{0,[J_n]})(\tilde\ga-\tilde\ga_ {[J_n]}) = o(1)$ as
%in \cite{cr13}.
%Lemma change of variable. Conclusion.

%s4.5 #&#
\subsection{Tools for density estimation}

The notation here follows the one introduced in Section~\ref{sec-dens-not},
in particular $\llVert\cdot\rrVert_L, W_n, R_n$ and $ {\mathcal{B}}$.

%Here write down lemma about upper bound on posterior behaviour of the
%Laplace transform of $t\psg f-f_0, \ga\psd_2$ under $\Pi^{D_n}$, for
%a general $\gamma$ (mind the centering at $f_0$ instead of $

%le3 #&#
\begin{lem} \label{lem-lapu}
Let $f_0$ belong to $ {\mathcal{F}}_0$. %=\cF_{m,M}$.
Let $\{a_n\}$ be a sequence of real numbers with $n a_n^2\ge1$, any
$n\ge1$.
Let $\{\Pi_n\}$ be a collection of priors on densities restricted to
the set
$\{f, h(f,f_0)\le a_n\}$.
Let $\{\gamma_n\}$ be an arbitrary sequence in $L^{\infty}[0,1]$. Set
$\tilde\gamma_n:= \gamma_n - P_{f_0}\gamma_n$. Suppose, for some
$m>0$ and %some $C_2>0$ and
all $n\ge1$,
\[
\llVert\tilde\gamma_n \rrVert_L \le m, \qquad\llVert
\tilde\gamma_n \rrVert_{\infty
} \le\bigl(4a_n
\log(n+1)\bigr)^{-1}.
\]
Then there exists $C>0$ depending on $m, \llVert f_0\rrVert_{\infty
}$ only %$N_0
such that for any $n\ge1$ and any $\llvert t\rrvert\le\log{n}$,
\[
E^{\Pi_n} \bigl[ e^{ t \sqrt{n}{\langle}f-f_0, \gamma_n {\rangle
}_2 } \given X^{(n)} \bigr] \le
e^{ Ct^2 + tW_n(\gamma_n)} \frac{ \int e^{\ell_n(f_t)-\ell_n(f_0) } \,d\Pi_n(f)
}{\int e^{\ell_n(f)-\ell_n(f_0)}\,d\Pi_n(f) },
\]
where $f_t$ is defined by $\log f_t = \log f - t\tilde\gamma_n/\sqrt
{n}- c(f e^{- t\tilde\gamma_n/\sqrt{n}})$.
\end{lem}

%%%%%%%% Proof of Laplace lemma %%%%%%%

\begin{pf}%{Proof of Lemma~\ref{lem-lapu}}
Denote $g=\log f$, $g_0=\log f_0$. From elementary algebra, it follows that
\begin{eqnarray*}
&& t\sqrt{n} {\langle}f-f_0, \gamma_n {
\rangle}_2 + \ell_n(f) - \ell_n(f_0)
\\
&&\qquad =  -\frac{n}2\biggl\llVert g - g_0 - \frac{t}{\sqrt{n}}
\tilde\gamma_n\biggr\rrVert_L^2 +
\sqrt{n}W_n\biggl(g - g_0 - \frac{t}{\sqrt{n}}\tilde
\gamma_n\biggr)
\\
&&\quad\qquad{}  + t\sqrt{n} {\mathcal{B}}(\tilde\gamma_n,f,f_0) +
R_n(f,f_0)+ \frac{t^2}{2}\llVert\tilde
\gamma_n\rrVert_L^2 + t W_n(
\tilde\gamma_n)
\\
&&\qquad =  \ell_n(f_t) - \ell_n(f_0)
+ \bigl[ t\sqrt{n} {\mathcal{B}}(\tilde\gamma_n,f,f_0)
+ R_n(f,f_0) - R_n(f_t,f_0)
\bigr]
\\
&&\quad\qquad{} + \frac{t^2}{2}\llVert\tilde\gamma_n\rrVert
_L^2 + t W_n(\tilde\gamma_n).
\end{eqnarray*}
Let us show that the bracketed term is small. From the definition of
$R_n, f_t$,
%E_{f_0}^n E^{\Pi^{D_n}}[e^{t\rn(f_{lk}-\hat{f}_{lk})} \given X]
%
\begin{eqnarray*}
&& R_n(f,f_0)-R_n(f_t,f_0)
\\
&&\qquad
= -\frac{t^2}{2}\llVert\tilde\gamma_n\rrVert
_L^2 + t\sqrt{n} {\langle}g-g_0,\tilde
\gamma_n{\rangle}_L + n \log F \bigl[e^{-t/\sqrt{n}\tilde\gamma_n}
\bigr].
\end{eqnarray*}
Next, expand the last logarithmic term. Using the assumption on $\llVert
\tilde\gamma_n\rrVert_\infty$ and $\llvert t\rrvert\le\log{n}$, the
absolute value
of the exponent in this term is at most $1/4$.
%First, since $f_0$ is bounded, $\left\Vert\tilde{\psi}_{lk}\right
%regime, this
%is at most $\sqrt{J_n}=o(\rn)$.
This enables us to expand successively the logarithm and exponential
functions, using the inequalities (the first is valid for $\llvert
x\rrvert\le1/4$),
\[
e^{-x} \le1-x + x^2,\qquad\log(1+x) \le x,
\]
leading to
\begin{eqnarray*}
\log F \bigl[e^{-t/\sqrt{n}\tilde\gamma_n} \bigr] & \le&\log F
\biggl[1-\frac{t}{\sqrt{n}}
\tilde\gamma_n + \frac
{t^2}{n}\tilde\gamma_n \biggr]
\\
& \le&\log\biggl[ 1- \frac{t}{\sqrt{n}}F\tilde\gamma_n +
\frac
{t^2}{n}F\tilde\gamma_n^2 \biggr]
\\
& \le&- \frac{t}{\sqrt{n}}F\tilde\gamma_n + \frac{t^2}{n}F_0
\tilde\gamma_n^2 + \frac{t^2}{n}(F-F_0)
\tilde\gamma_n^2.
\end{eqnarray*}
The last term, using Lemma~\ref{lem-tech} and
$h(f,f_0) \llVert\tilde\gamma_n\rrVert_{\infty} \le1$ together
with $\llVert\tilde\gamma_n\rrVert_2^2 \lesssim\llVert\tilde\gamma
_n\rrVert
_L^2\lesssim m$, is a $O(t^2/n)$.
%For $l\le J_n$ and $\veps_n=n^{-\alpha/(2\alpha+1)}$ with $
%C$ for some finite constant $C$. This implies that the last term is
%$O(1)t^2/n$.
On the other hand,
\[
F\tilde\gamma_n=(F-F_0)\tilde\gamma_n = {
\biggl\langle}\frac
{f-f_0}{f_0}, \tilde\gamma_n {\biggr
\rangle}_L = {\langle}g-g_0, \tilde\gamma_n
{\rangle}_L + {\mathcal{B}}( \tilde\gamma_n,f,f_0).
\]
Combine the previous results to obtain the desired bound.
\end{pf}

%Here write down concentration rate for the density estimation in the
%squared $L^2$-norm (tools from the density paper) and a preliminary
%sup-norm rate.

%le4 #&#
\begin{lem} \label{lem-inter}
Consider the log-density prior (\ref{def-prior-log}). Suppose
$g_0=\log f_0$ belongs to $ {\mathcal{C}}^\alpha$, with $\alpha
>1/2$. Suppose (\ref{rate-hell}) holds and let $\varepsilon_n, \zeta
_n$ be defined as below (\ref{rate-hell}).
Then $M$ large enough,
\[
E_{f_0}^n\Pi\bigl[f\dvtx  \llVert f-f_0\rrVert
_2\le M\varepsilon_n, \llVert f-f_0\rrVert
_{\infty} \le M\zeta_n \given X^{(n)}\bigr] \to1.
\]
\end{lem}
\begin{pf}
Obtaining this result could be done following the arguments in \cite{gn11}.
Here we use instead an approach from \cite{rr12}, extending their
argument on the second moment of $\log(f/f_0)$ to further get a rate
$\zeta_n$ for $\llVert\log(f/f_0)\rrVert_\infty$.

Since (\ref{rate-hell}) is assumed, one can restrict to the event $\{
f\dvtx  h(f,f_0)\le\varepsilon_n\}\subset\{f\dvtx  \llVert f-f_0\rrVert
_2\lesssim
\varepsilon_n\}$.
We have (see, e.g., Lemma 8 in \cite{gvv07}),
\[
\bigl\llVert\log(f/f_0) \bigr\rrVert_{2}^2
\lesssim h^2(f,f_0) \bigl( 1+\log\llVert
f/f_0\rrVert_{\infty
} \bigr).
\]
The last term is bounded by a constant times $\varepsilon_n^2 (1+\log
\llVert T-c(T)-g_0 \rrVert_\infty)$ by assumption, where $g_0=\log
f_0$. Next,
one writes
%one bounds this last quantity using its explicit expression
%
\begin{eqnarray*}
&& \bigl\llVert T-c(T)-g_0 \bigr\rrVert_\infty
\\
&&\qquad = \biggl
\llVert\sum_{l,k} {\bigl\langle}T-c(T)-g_0,\psi_{lk}{\bigr\rangle}_2
\psi_{lk} \biggr
\rrVert_{\infty}
\\
&&\qquad \le\sum_{l\le L_n} 2^{l/2} \max
_{k} \bigl\llvert{\bigl\langle}T-c(T)-g_0,
\psi_{lk}{\bigr\rangle}_2 \bigr\rrvert+ \sum
_{l > L_n} 2^{l/2} \max_{k}\bigl
\llvert{\langle}g_0, \psi_{lk} {\rangle}_2
\bigr\rrvert.
\end{eqnarray*}
Since $g_0$ is $\alpha$-H\"{o}lder,\vspace*{1pt} the last term is of the order
$\varepsilon_{n,\alpha}^*=o(1)$. For the middle term, Cauchy--Schwarz
inequality yields the bound
$2^{L_n/2} \llVert T-c(T)-g_0\rrVert_2$, using $\sum_{l\le L_n} 2^l
\lesssim
2^{L_n}$ and bounding
the maximum of squares by the sum. Deduce that
\[
\bigl\llVert\log(f/f_0) \bigr\rrVert_{2}^2
\lesssim\varepsilon_n^2 + \varepsilon_n^2
2^{L_n/2} \bigl\llVert\log(f/f_0) \bigr\rrVert
_{2} \lesssim\varepsilon_n^2 +
\varepsilon_n^2 2^{L_n} \bigl\llVert
\log(f/f_0) \bigr\rrVert_{2}^2.
\]
Since $\alpha>1/2$, one has $\varepsilon_n^2 2^{L_n}=o(1)$. So
gathering the $L^2$-norm terms on the same side of the inequality one
obtains $\llVert\log(f/f_0) \rrVert_{2}^2 \lesssim\varepsilon_n^2$.
%% on the
%event $\{f, h(f,f_0)\le\veps_n\}$.
Also, along the way
%under the event $\{f, h(f,f_0)\le\veps_n\}$,
we have obtained the bound
\[
\bigl\llVert T-c(T)-g_0 \bigr\rrVert_\infty
\lesssim2^{L_n/2} \bigl\llVert T-c(T)-g_0 \bigr\rrVert
_2 + \varepsilon_{n,\alpha}^* \lesssim2^{L_n/2}
\varepsilon_n = \zeta_n.
\]
Now the squared $L^2$-norm of $f-f_0$ can be expressed as
\[
\int_0^1 (f-f_0)^2
= \int_0^1 f_0^2
\bigl(e^{T-c(T)-g_0} - 1\bigr)^2.
\]
The inequality $\llvert e^x -1 \rrvert\le C\llvert x\rrvert$, valid
for $x$ in a compact subset
of $\mathbb{R}$ and $C$ a large enough constant, implies %again in
%case $f$ is in an $\veps_n$-Hellinger neighborhood of $f_0$,
%
\[
\int_0^1 (f-f_0)^2
\le C^2 \int_0^1
f_0^2\bigl(T-c(T)-g_0\bigr)^2
\lesssim\bigl\llVert T-c(T)-g_0\bigr\rrVert_2^2
\lesssim\varepsilon_n^2.
\]
Similarly, since $\llVert f_0\rrVert_\infty<\infty$, one obtains
$\llVert f-f_0\rrVert
_\infty\lesssim\zeta_n$. %on $\{f, h(f,f_0)\le\veps_n\}$.
\end{pf}

%s4.6 #&#
\subsection{Other lemmas}

Given $R>0$, let $B^{\alpha}_{\infty,\infty}(R)$ denote the centered
ball of $B^{\alpha}_{\infty,\infty}[0,1]$ of radius $R$ for the norm
$\llVert\cdot\rrVert_{\infty,\infty,\alpha}$ given in Section~\ref{wave}.

%le5 #&#
\begin{lem} \label{lem-hoelder}
There exists a constant $C>0$ such that for any $f,g$ in $B^{\alpha
}_{\infty,\infty}(R_f)$ and $B^{\alpha}_{\infty,\infty}(R_g)$,
respectively, the product $fg$ belongs to\break  $B^{\alpha}_{\infty,\infty
}(CR_fR_g)$. If $f$ belongs to $ {\mathcal{C}}^\alpha[0,1]$ and is
bounded away from $0$ then $f^{-1}$ belongs to $ {\mathcal
{C}}^{\alpha}[0,1]$. Moreover, for any indexes $l,k$,
\[
\llVert\psi_{lk}\rrVert_{\infty,\infty,\alpha} = 2^{l(1/2 +\alpha
)}.
\]
\end{lem}
\begin{pf}
The first claim follows from the main result of Section~2.8.3 in \cite
{t83} (strictly speaking the last result is for functions on $\mathbb
{R}$, %whereas we consider functions on $[0,1]$,
but the latter functions can be shown to be restrictions to $[0,1]$ of
elements of $B_{\infty,\infty}^\alpha$ whose norm is equivalent to
the one of the restriction; see \cite{n07} Proposition~2 for a similar
argument for Sobolev spaces). The second claim is a simple computation
using the definition of H\"{o}lder spaces. For the last claim, one uses
the characterisation of $B^{\alpha}_{\infty,\infty}$ in terms of
wavelet coefficients from Section~\ref{wave}, which yields $\llVert
\psi
_{lk}\rrVert_{\infty,\infty,\alpha} =
\max_{l',k'} 2^{l'(1/2 +\alpha)} {\langle}\psi_{lk}, \psi
_{l'k'} {\rangle}_2 = 2^{l(1/2 +\alpha)}$.
\end{pf}

%le6 #&#
\begin{lem} \label{lem-tech}
Let $f,f_0$ be two densities on $[0,1]$ such that $f_0$ is bounded.
For any $g\in L^2[0,1]$ such that $h(f,f_0)\llVert g\rrVert_{\infty}
\le C_1$ and
$\llVert g\rrVert_2 \le C_2$, for some constants $C_1, C_2>0$,
\[
\bigl\llvert(F-F_0)g^2\bigr\rrvert\le
C_1^2+C_1\sqrt{4C_2
\llVert f_0\rrVert_{\infty}+ C_1^2}.
\]
\end{lem}
\begin{pf}
Denote $\Sigma:= \int_0^1 \llvert f-f_0\rrvert g^2 $. Then by
Cauchy--Schwarz inequality,
\begin{eqnarray*}
\Sigma^2 & \le&2 h(f,f_0)^2 \int
_0^1(f+f_0) g^4
\\
& \le&2 h(f,f_0)^2 \bigl[ \llVert g\rrVert
_{\infty}^2 \Sigma+ 2 \llVert f_0\rrVert
_{\infty} \llVert g\rrVert_{\infty}^2\llVert g\rrVert
_{2}^2 \bigr]
\\
& \le&2C_1^2\bigl[\Sigma+ 2 C_2\llVert
f_0\rrVert_{\infty}\bigr].
\end{eqnarray*}
This implies that $\Sigma$ is less than the largest root of the
polynomial $X^2-2C_1X+4C_1^2\llVert f_0\rrVert_{\infty}$, which can
be expressed
in terms of
$C_1, \llVert f_0\rrVert_{\infty}$. %Standard inequalities yield the
%desired
%bound.
\end{pf}

%le7 #&#
\begin{lem} \label{lem-Gamma}
Let $f_0\in{\mathcal{F}}_0$ and $g_0=\log f_0$, and let $\Gamma
^{L_n}$ be defined by (\ref{sto}), with $ {\mathcal{A}}_{l,k}$ any
elements of $L^\infty[0,1]$ such that\vspace*{1pt} there exists constants $c_1,
c_2$ with, for any $l,k$ with $k < 2^l \le2^{L_n}$, any $n\ge2$,
\[
\llVert{\mathcal{A}}_{l,k} \rrVert_{\infty} \le
c_1\sqrt{n/\log{n}},\qquad\llVert{\mathcal{A}}_{l,k}
\rrVert_2 \le c_2.
\]
Then for any $n\ge2$ and $L_n$ defined by (\ref{def-jn}), it holds
\[
E_{f_0}^n \bigl\llVert\Gamma^{L_n} -
g_0^{L_n}\bigr\rrVert_{\infty} \lesssim\varepsilon
_{n,\alpha}^*.
\]
\end{lem}
\begin{pf}
We proceed exactly as for the proof of the maximal inequality in
Section~\ref{firstex}.
For any $t>0$,
\[
E_{f_0}^n \bigl\llVert\Gamma^{L_n} -
g_0^{L_n}\bigr\rrVert_{\infty} \le\frac{1}{\sqrt{n}}
\sum_{l\le L_n} \frac{2^{{l}/{2}}}{t} \log\sum
_{k=0}^{2^l -1} E_{f_0}^n\bigl[
e^{t W_n( {\mathcal
{A}}_{l,k})} + e^{-t W_n( {\mathcal{A}}_{l,k})} \bigr].
\]
We have $W_n( {\mathcal{A}}_{l,k})=\mathbb{G}_n( {\mathcal
{A}}_{l,k})$ and bounds on exponential moments of the last
empirical quantity are well known. From Laplace transform controls, one
gets, for any real $s$,
\[
E_{f_0}^n\bigl[ e^{s W_n( {\mathcal{A}}_{l,k})} \bigr] \le
e^{ (s^2/2) [\int{\mathcal{A}}_{l,k}^2 f_0 ]
e^{\llvert s\rrvert\llVert{\mathcal{A}}_{l,k} \rrVert_{\infty
}/\sqrt{n}} }.
\]
Let us choose $t=\sqrt{l}\lesssim\sqrt{L_n}$. Under the conditions
of the lemma the last display with $s=t$ or $s=-t$ is bounded above by
$e^{Ct^2}$. This leads to the bound $E_{f_0}^n \llVert\Gamma^{L_n} -
g_0^{L_n}\rrVert_{\infty} \lesssim\varepsilon_{n,\alpha}^*$.
\end{pf}

%le8 #&#
\begin{lem} \label{lem-exa}
Let $\varphi=\varphi_G$ and $\sigma_{l}=2^{-l(1/2 +\gamma)}$
for all $l\le L_n$, and any given $0< \gamma\le\alpha-1/4$. Then
(\ref{cond-sil})--(\ref{rate-hell}) hold. The same applies for
$\varphi=\varphi_{H,\tau}$ and $\sigma_{l}=2^{-l\alpha}$ for all
$l\le L_n$,
for any given value of the parameter $0\le\tau<1$.
\end{lem}
\begin{pf}
The first result is a minor adaptation of Theorem~4.5 in \cite{vvvz},
where the authors consider a cut-off at $n^{1/(2\alpha+1)}$ instead of
$2^{L_n}=h_n^{-1}=(n/\log n)^{1/(2\alpha+1)}$ (equality up to fixed
multiplicative constants). Taking the cut-off at $2^{L_n}$ only changes
logarithmic factors in their argument. In the log-Lipschitz case, one
adapts Theorem~2.1 in \cite{rr12}. There
are two points to note. First, taking $2^{L_n}$ instead of
$n^{1/(2\alpha+1)}$ induces only, again, an extra logarithmic power in
the rate. Second, strictly speaking the authors in \cite{rrr12}
consider wavelets on the interval via periodisation, which imposes
conditions at the boundary (periodic Besov spaces), conditions which
can be dropped when using the CDV wavelet basis. % enables one to
%consider the spaces without the boundary conditions, without any other
%change to the arguments.
Explicit (re)derivation of the previous two results is omitted.
\end{pf}

%le9 #&#
\begin{lem} \label{lem-max}
Let $\varphi=\varphi_G$ and $\sigma_{l}$ satisfy (\ref{cond-sil}).
Then the prior $\Pi$
defined by (\ref{def-prior-log}), with $L_n$ as in (\ref{def-jn}) and
$\varepsilon_n$ as below (\ref{rate-hell}), satisfies, for $C>0$
large enough and
any fixed given integer $K$,
\[
E_{f_0}^n\Pi\Bigl[ \max_{\lambda\le K,\mu} \bigl
\llvert{\langle}T,\psi_{\lambda\mu} {\rangle}\bigr\rrvert\le C\sqrt{n}
\varepsilon_n\big\given X^{(n)} \Bigr] \to1.
\]
\end{lem}
\begin{pf}
The\vspace*{1pt} maximum in the display of the lemma only involves a finite number
of terms $\lambda,\mu$ with $\lambda\le K$, $\mu\le2^{K}-1$, and
these terms are Gaussian with variances $\sigma_{l}^2$ bounded above
by positive constants. Thus, by Gaussian concentration one gets
\[
\Pi\Bigl[ \max_{\lambda\le K,\mu} \bigl\llvert{\langle}T,
\psi_{\lambda\mu
} {\rangle}\bigr\rrvert> C\sqrt{n}\varepsilon_n
\Bigr] \le e^{-cn\varepsilon
_n^2},
\]
where $c$ can be made arbitrarily large by taking $C$ large enough.
Next, one applies Lemma 1 in \cite{gvni}. To do so, one needs to bound
from below the prior probability of a Kullback--Leibler neighborhood of
$f_0$ of size $\varepsilon_n$ by $e^{-dn\varepsilon_n^2}$ for some
$d>0$. This follows from the conclusion of Theorem 5 in \cite{vvvz},
which (modulo the fact that our $\varepsilon_n$ is within a
logarithmic factor of theirs, as noted in Lemma~\ref{lem-exa} to
accommodate our slightly different choice
of cut-off $2^{L_n}$), provides the bound $\Pi[\llVert g-g_0\rrVert
_{\infty
}<4\varepsilon_n)\ge e^{-n\varepsilon_n^2}$. Switching from the
sup-norm on $g-g_0$ to the Kullback--Leibler divergence between $g$ and
$g_0$ follows from Lemma 3.1
in \cite{vvvz}.
\end{pf}

%le10 #&#
\begin{lem} \label{lem-histr}
Suppose $f_0$ belongs to $ {\mathcal{F}}_0\cap{\mathcal
{C}}^\alpha[0,1]$, for $0<\alpha\le1$.
Let $\Pi$ be a prior on histogram densities defined by (\ref{prk}).
Then, for $M$ large enough
\[
E_{f_0}^n\Pi\bigl[f, h(f,f_0)\le M(
\log{n}/n)^{-\alpha/(2\alpha +1)} \given X^{(n)} \bigr] \to1.
\]
\end{lem}
\begin{pf}
One can proceed as in the proof of Theorem 2 in \cite{s07}. It suffices
to use as sieve the set of dyadic density histograms $ {\mathcal
{H}}_{L_n}$, so that \mbox{$\Pi[ {\mathcal{H}}_{L_n}^c]=0$}. Entropy and
prior mass conditions can then be verified with the same rate $(\log
n/n)^{\alpha/(2\alpha+1)}$;
see also the remark before Theorem 2.4 in \cite{ghosal01}.
%The rate follows from an application of the general rate theorem in
%6.1 on Dirichlet weights. If the weight sequence
\end{pf}

%le11 #&#
\begin{lem} \label{lem-jaco}
Let $\Delta(\zeta)$ be the Jacobian of the change of variables
$\omega\to\zeta$ given by (\ref{cvarsim}) over the unit simplex
$ {\mathcal{S}}_{2^{L_n}}$. It holds, with $f(\zeta)=
2^L \sum_{\mu=0}^{2^L-1} \zeta_\mu\mathbh{1}_{I_k^L}$,
\[
\Delta(\zeta) = \prod_{\mu=0}^{2^{L_n}-1}
\frac{e^{t \gamma
_{n,\mu}/\sqrt{n}}}{\int_0^1 e^{t\gamma_n(x)/\sqrt{n}}f(\zeta
)(x)\,dx}.
\]
\end{lem}
\begin{pf}
This follows from elementary calulations; see \cite{cr13}.
\end{pf}

\section*{Acknowledgements}
This work started from ideas developed in \cite
{CN12} in common with Richard Nickl. His numerous helpful comments are
warmly acknowledged. I am also indebted to Judith Rousseau for early
explanations on \mbox{semi}parametric analysis of histograms priors. Although
the treatment for histograms here is a bit different, it owes to her
ideas, and we refer to \cite{cr13} for a systematic treatment of
histograms for fixed functionals. I would like to thank Catia
Scricciolo for early discussions and G\'erard Kerkyacharian and
Dominique Picard for helpful comments, as well as the referees and the
Associate Editor for insightful and constructive comments.

%suskaldyti doi

% imsref loaded by linak, 2014-07-28 12:26:02
%
% imsref loaded by linak, 2014-08-12 10:32:04

\printaddresses
\end{document}